\documentclass[reqno,11pt]{amsart}

\newtheorem{theorem}{Theorem}[section]
\newtheorem{lemma}[theorem]{Lemma}
\newtheorem{corollary}[theorem]{Corollary}

\theoremstyle{definition}
\newtheorem{assumption}[theorem]{Assumption}

\theoremstyle{remark}
\newtheorem{remark}[theorem]{Remark}
\newtheorem{example}[theorem]{Example}

\makeatletter 
 \def\dashint{%
 \operatorname%
 {\,\,\text{\bf--}\kern-.98em\DOTSI\intop\ilimits@\!\!}}
\def\dashnorm{\,\,\text{\bf--}\kern-.5em\|}
\def\ninf{\qopname\relax\@empty{inf\phantom{p}\!\!\!}}

 \makeatother

\def\sft{{\sf t}}

\newcommand\bB{\mathbb{B}}
\newcommand\bC{\mathbb{C}}

\newcommand\bR{\mathbb{R}}
\newcommand\bS{\mathbb{S}}

\newcommand \bW{\mathbb{W}}

\newcommand\cF{\mathcal{F}}

\newcommand \cL{\mathcal{L}}

\newcommand\frN{\mathfrak{N}}

\newcommand{\loc}{{\rm loc}\,}

 \newcommand{\mysection}[1]{\section{#1}
 \setcounter{equation}{0}}
 
\newcommand{\esssup}{\operatornamewithlimits{ {ess\,sup}}}

\begin{document}

\title[On weak  and strong solutions]
{On weak  and strong solutions of time inhomogeneous
It\^o's equations with VMO diffusion and Morrey drift}
\author{N.V. Krylov}

\email{nkrylov@umn.edu}
\address{School of Mathematics, University of Minnesota, Minneapolis, MN, 55455}
 
\keywords{Strong solutions, time inhomogeneous equations, Morrey drift}
 
\subjclass{60H10, 60J60}

\begin{abstract} 
We prove the existence of weak solutions
of It\^o's stochastic time dependent equations with
irregular diffusion and drift terms
of Morrey  spaces. Weak uniqueness (generally conditional)
and 
 a conjecture  pertaining to strong solutions are also
discussed.  Our results are new even
if the drift term vanishes. 
\end{abstract}

\maketitle

\mysection{Introduction}
                                                  \label{section 3.11.1}

This paper is a natural continuation of
the second part of \cite{7}, where the weak uniqueness is proved for the time homogeneous stochastic equations with Morrey drift. Here we do the same
in the time inhomogeneous case  
and also explore some features of strong solutions.

Let $\bR^{d}$ be a $d-$dimensional Euclidean space of points
$x=(x^{1},...,x^{d})$ with $d\geq2$. Let $(\Omega,\cF,P)$ be a 
complete probability space,
carrying a $d $-dimensional Wiener process
  $w_{t}$. Fix $\delta\in(0,1]$ and denote by $\bS_{\delta}$ the set of $d\times d$ symmetric matrices
whose eigenvalues lie in $[\delta,\delta^{-1}]$.

Assume that on $\bR^{d+1}=\{(t,x):t\in\bR,x\in\bR^{d}\}$ we are given Borel $\bR^{d}$-valued 
function  $b=(b^{i})$ and $\bS_{\delta}$-valued $\sigma 
=(\sigma^{ij})$.
We are going
to   investigate
the equation
\begin{equation}
                                         \label{6.15.2}
 x_{s}=x +\int_{0}^{s}\sigma (t+u,x_{u})\,dw _{u}+
\int_{0}^{s}b(t+u,x_{u})\,du.
\end{equation}

We are interested in the so-called weak solutions, that is
solutions that are not necessarily $\cF^{w}_{s}$-measurable,
where $\cF^{w}_{s}$ is the completion of $\sigma(w_{u}:u\leq s)$. We present
sufficient conditions for the equation to have such solutions on appropriate probability spaces and investigate uniqueness of their distributions.

After the classical work by K. It\^o showing that there exists
a unique (strong) solution of \eqref{6.15.2} if $\sigma$ and $b$
are Lipschitz continuous in $x$ (may also depend on   $\omega$ and the nondegeneracy of $\sigma$ is not required),
  much effort was  applied to relax these Lipschitz conditions.   The first author who achieved a considerable progress was A.V. Skorokhod
\cite{Sk_61} who proved the solvability assuming
only the continuity of $\sigma$ and $b$ 
in $x$ (which may depend
on $t$ and again without nondegeneracy).
Then by using the Skorokhod method and Aleksandrov
estimates the author proved in \cite{Kr_69}
and \cite{Kr_77} the solvability
for the case of {\em just measurable $\sigma$ and  
bounded $b$}  under the nondegeneracy assumption.
Stroock and Varanhan \cite{SV_79} among other things
not only proved
the solvability for the coefficients uniformly 
continuous in $x$, but also proved the uniqueness
of their distributions.

A few words about strong solutions.
In   case $d =1$ T. Yamada and S. Watanabe \cite{YW_71} relaxed
the Lipschitz condition on $\sigma$ to the H\"older $(1/2)$-condition
(and even slightly weaker condition) and kept $b$ Lipschitz
(slightly less restrictive). Much attention was paid to equations
with continuous coefficients
satisfying the so-called monotonicity conditions
(see, for instance, \cite{Kr_84} and the references therein).

T. Yamada and S. Watanabe \cite{YW_71} also put forward
a very strong theorem, basically, saying that 
pathwise uniqueness and existence of weak solutions
imply the existence of strong
solutions. After that
the majority of papers on the subject  
are using their theorem.
 S. Nakao (\cite{Na_72}) proved
the strong solvability in time homogeneous case
 if $d=1$ and $\sigma$ is bounded away from zero and infinity
 and is
locally of bounded variation. He also assumed that $b$ is bounded,
but from his arguments it is clear that the summability of $|b|$
suffices.

A. Veretennikov was the first author who in \cite{Ve_80} 
not only proved the existence of strong solutions
in the time inhomogeneous multidimensional case  when $b$ is bounded,
but also considered the case of $\sigma^{k}$ in Sobolev class,
namely,   $\sigma^{k}_{x}\in L_{2d,\loc}$. He used A. Zvonkin's method (see
\cite{Zv_74}) of transforming the equation in such a way that the
drift term disappears.
In \cite{XXZZ_20} (also see the references
there) the result of Veretennikov is
extended to the case of $\sigma$ uniformly continuous in $x$ and $\sigma_{x},b\in L_{p,q}$
with, perhaps, different $p,q$ for $\sigma_{x}$ and $b$ satisfying
\begin{equation}
                           \label{6.12.1}
\frac{d}{p}+\frac{2}{q}< 1.
\end{equation}
In that case much information is available,
we refer the reader
to \cite{Zh_11},  \cite{Zh_20}, \cite{XXZZ_20}, and the references therein.  

Even the case when   $\sigma$ is constant and the process is nondegenerate
attracted
very much attention.   
M. R\"ockner and the author in \cite{KR_05} proved,
among other things, the existence of strong solutions 
when $b\in L_{p,q}$ under condition \eqref{6.12.1}.
If $b$ is bounded  A. Shaposhnikov (\cite{Sh_14}, \cite{Sh_17})
proved the so called path-by-path uniqueness, which, basically,
means that for almost any trajectory $w_{t}$ there is only one
solution (adapted or not). This result was already announced
by A. Davie before with a very entangled proof which left many doubtful.

In a fundamental work by L.~Beck, F.~Flandoli, M.~Gubinelli, and M.~Maurelli
(\cite{BFGM_19}) the authors investigate such equations from
 the points of view of It\^o stochastic equations, stochastic transport 
equations, and stochastic continuity equations. Their article
contains an enormous amount of information and a vast references list. 
In what concerns our situation they require ($\sigma$ constant
and the process is nondegenerate) what they call LPS-condition (slightly imprecise): $b\in L_{p,q,\loc}$,
$q<\infty$, with equality in \eqref{6.12.1} in place of $<$,
or $p=d$ and $\|b\|_{L_{p,\infty}}$ to be sufficiently small, or else that $b(t,\cdot)$
be continuous as an $L_{d}(\bR^{d})$-function,
and they prove strong solvability and strong uniqueness
(actually, path-by-path-uniqueness which is stronger) only for {\em almost
all\/} starting points $x$. 
Concerning the strong solutions
starting from any point $x$
in the {\em time dependent\/} case
with singular $b$  and constant $\sigma$ probably the best results belong to R\"ockner and Zhao \cite{RZ_21}, where, among very many other things, they prove existence and uniqueness of strong solutions of
equations like \eqref{6.15.2} with $b\in L_{p,q}$ and $p,q<\infty$, 
with equality in \eqref{6.12.1} in place of $<$, or when $b(t,\cdot)$
is continuous as an $L_{d}(\bR^{d})$-function.
  In the recent paper \cite{KM_23}
the authors consider the equations with constant
$\sigma$ and the drift term more general than in
\cite{RZ_21} in some respects, however, not completely
covering all possible drifts in \cite{RZ_21}.
 In Section \ref{section 2.28.1}  of the current paper
 we present a conjecture and its intuitive justification
concerning the existence of strong solutions
in the case of {\em nonconstant\/}
 diffusion and singular drift which is valid
 for the results in \cite{RZ_21}. 

Coming back to weak solutions it is worth saying
that
restricting
the situation to the one when $\sigma$ and $b$
are independent of time allows one to
relax the above conditions significantly
further, see, for instance, \cite{KS_19}
and the 
references therein. 

This paper
is close to  \cite{Ki_23}, where the drift
term is   more general  
in some respects, however, not completely
covering the results in \cite{RZ_20} or our results
in case $b^{M}\equiv0$, 
 when $\sigma$
is the unit matrix and the weak uniqueness
is proved among the solutions which are constructed
by approximating the initial coefficients by smooth ones. Our uniqueness theorem and uniqueness
theorems in \cite{RZ_20} are also conditional.
We prove uniqueness only in the class of solutions
({\em which is proved to be nonempty\/}) admitting certain
estimates,   however,
in Remark \ref{remark 10.30.10} we 
mention a sufficient analytic
condition on $b$ when the unconditional weak uniqueness
holds. 

In Remarks 
\ref{remark 1.5.1}  and 
\ref{remark 10.30.10}  we compare our results with some
of those in an excellent paper by
R\"ockner and Zhao \cite{RZ_20} and 
 refer the reader to \cite{BFGM_19}, \cite{Ki_23}, \cite{RZ_20}   for
 very good reviews  of the recent history 
 of the problem.
 
 By the way,
G. Zhao (\cite{Zh_20_1}) gave an example showing
that, if in condition \eqref{6.12.3}   we replace $\rho^{-1}$ 
with $\rho^{-\alpha}$, $\alpha>1$,   weak uniqueness
may fail even in the time homogeneous case
and unit diffusion.  According to 
\cite{MG_23},   assuming that $b\in L_{p,q}$ with $d/p +1/q \leq 1$
alone does not guarantee weak uniqueness even with unit diffusion
(the existence is known).
In Example 2.1 of \cite{Kr_20}, for any $\varepsilon>0$ we have
$b\in L_{p,q}$ with $d/p +1/q \geq 1+\varepsilon$ and there are no solution
of \eqref{6.15.2} with unit diffusion and $(t,x)=0$ at all.
In Example 3 of \cite{KZ_75} it is given an equation
$dx_{t}=\sigma(x_{t})\,dw_{t}$ in $d=2$ with $\sigma\sigma^{*}
=(\delta^{ij})$ such that it has unique and strong solutions
for any starting point apart from the origin. If the
starting point is the origin, only weak solutions exist.
All these examples show that we are dealing with quite
delicate problems, many of which are to date far
from being settled in the most satisfactory way.

Here is an example in which we prove existence 
(and conditional uniqueness) of weak solutions:
$\sigma=2(\delta^{ij})+I_{x \ne0}\zeta(x )\sin(\ln|\ln |x |)$
 (quite discontinuous),
where $\zeta$ is any smooth symmetric $d\times d$-matrix valued 
function vanishing for $|x|>1/2$ and
satisfying $|\zeta|\leq 1$, and $|b|=\gamma/|x|$
with $\gamma$ sufficiently small.
Another example of $b$ is when
$|b|= \gamma(|x|+\sqrt{|t|})^{-1}I_{|x|<1,|t|<1}$
with $\gamma$ sufficiently small.
Both examples of $b$ are admissible in
\cite{Ki_23} and inadmissible in 
  \cite{RZ_20}. In both articles $\sigma$ is constant.
   By the way, it is well known that
  the equation $dx_{t}=
  dw_{t}-|x_{t}|^{-1}b(x_{t})\,dt$, where $b(x)
  =(d/2)x/|x|$, with initial data $x_{0}=0$
  does not have solutions, so that in the above examples
  $\gamma$ indeed should be sufficiently small.

 The paper is organized as follows.
In Section \ref{section 12.20.1} we set up
the problem and formulate our assumptions.
In Section \ref{section 3.4.1} we prove the solvability
of \eqref{6.15.2} and prove It\^o's formula
for the solutions. This formula is the main tool
to prove in Section \ref{section 3.4.2}
the weak uniqueness of our solutions and the fact
that they
constitute a strong Markov process with strong Feller
resolvent. In Section \ref{section 2.28.1}
we prove a uniqueness result for  strong solutions
  and, in case
of smooth $\sigma$ and $b$,  prove a formula 
which allows us to state a conjecture
about the existence of strong solutions
with $D\sigma$ almost in VMO and $b$ in a Morrey class.
Proving this conjecture is our next quite
challenging project. We intend to use some ideas
from \cite{Kr_21} and \cite{1}. However,
as the authors of \cite{RZ_21} write
``It may be also possible to follow the same procedure as 
in \cite{Kr_21} to study
the time-inhomogeneous case, but 
one encounters a lot of difficulties 
due to the fact that there is no
good PDE theory for equations with 
such kind of first order terms so far.''
Indeed, in the past a few attempts of the author  
to prove the conjecture failed.

We conclude the introduction by some notation.
We set  
$$
D_{x^{i}}u=D_{i}u=u_{x^{i}}=\frac{\partial}{\partial x^{i}}u,\quad Du=u_{x}=(D_{i}u),
$$
$$
  D_{ij}u=u_{x^{i}x^{j}}=D_{i}D_{j}u ,\quad u_{xx}=(D_{ij}u),
$$
$$
u_{x^{i}\eta^{j}}=D_{x^{i}\eta^{j}}u=
D_{x^{i}}D_{\eta^{j}}u,\quad u_{x\eta}=(D_{x^{i}\eta^{j}}u),
$$
$$
u_{\eta^{i}\eta^{j}}=D_{\eta^{i}\eta^{j}}u=
D_{\eta^{i}}D_{\eta^{j}}u,\quad
u_{\eta\eta}=(D_{\eta^{i}\eta^{j}}u),
$$
$$
 \partial_{t}u=\frac{\partial}{\partial t}u,\quad
u_{(\eta)}=\eta^{i}u_{x^{i}}.
$$
 If $\sigma =(\sigma^{ij...})$ by $|\sigma|^{2}$
we mean the sum of squares of all entries.

For $p\in[1,\infty)$, 
and domain $\Gamma\subset \bR^{d}$ by   $L_{p}(\Gamma)$  we mean the space 
of Borel (real-, vector-, matrix-valued...)
 functions on   $\Gamma $  with finite norm given by
$$
 \|f\|_{L_{p}(\Gamma)}^{p}=\int_{\Gamma}|f(x)|^{p}\,dx .
$$
Set $L_{p}=L_{p}(\bR^{d})$.

For $p,q\in[1,\infty)$ and domain $Q\subset\bR^{d+1}$ by $L_{p,q}(Q)$
we mean the space of Borel (real-, vector- or matrix-valued)
 functions on $Q$   with finite norm given by
$$
\|f\|_{L_{p,q}(Q)}^{q}=\|fI_{Q}\|_{L_{p,q}}^{q}
=\int_{\bR}\Big(\int_{\bR^{d}}|fI_{Q}(t,x)|^{p}\,
dx\Big)^{q/p}\,dt.
$$
Set $L_{p,q}=L_{p,q}(\bR^{d+1})$.
By $W^{2}_{p}$ we mean the space 
of Borel functions $u$ on $\bR^{d}$ whose Sobolev derivatives 
$u_{x}$ and $u_{xx}$ exist and $u,u_{x},u_{xx}\in L_{p}$.
The  norm in $W^{2}_{p}$ is given by
$$
\|u\|_{W^{2}_{p}}=\|u_{xx}\|_{L_{p}}+\|u \|_{L_{p}}.
$$ 
Similarly  $W^{1}_{p}$ is defined. As usual, we write $f\in L_{p,\loc}$
if $f\zeta\in L_{p}$ for any $\zeta\in C^{\infty}_{0}$ ($=C^{\infty}_{0}
(\bR^{d})$).  

By $W^{1,2}_{p,q}(Q)$ we mean the collection
of $u$ such that $\partial_{t}u$, $u_{xx}$, $u_{x}$, $u
\in L_{p,q}(Q)$. The norm in $W^{1,2}_{p}(Q)$
is introduced in an obvious way.
We abbreviate $W^{1,2}_{p,q } =W^{1,2}_{p,q}(\bR^{d+1})$.

If a Borel $\Gamma\subset \bR^{d}$, by $|\Gamma|$ we mean its Lebesgue
measure,
$$
\dashint_{\Gamma}f(x)\,dx=\frac{1}{|\Gamma|}
\int_{\Gamma}f(x)\,dx.
$$
Similar notation is used for $\Gamma\subset
\bR^{d+1}$.

  Introduce
$$
B_{R}(x)=\{y\in\bR^{d}:|x-y|<R\},\quad B_{R} =B_{R}(0)
$$
and let $\bB_{R}$ be the collection of balls of radius $R$. Also let  
$$
C_{\tau,\rho}(t,x)=[t,t+\tau)\times B_{\rho}(x),\quad C_{\rho}...=C_{\rho^{2},\rho}...,\quad C_{\rho}=C_{\rho}(0,0),
$$
and let $\bC_{\rho}$ be the collection of
$C_{\rho}(t,x)$.

In the proofs of our results
we use various (finite) constants called $N$ which
may change from one occurrence to another
and depend on the data only in the same way as
  indicated in the statements
of the results.

\mysection{Setting and assumptions}
 
                 \label{section 12.20.1}

  Set $a=\sigma^{2}$ and fix some $ r_{a},
 r_{b} \in(0,\infty)$.
The values of $\theta, \hat b_{M} >0$
in Assumption \ref{assumption 6.3.1} below,
that is supposed to hold throughout this section, will be specified later.  
\begin{assumption}
                     \label{assumption 6.3.1}

(i) For $\rho\leq r_{a}$  
\begin{equation}
                                 \label{6.3.01}
a^{\# }_{\rho}:= \sup_{C\in\bC_{\rho}}\dashint_{C}|a (t,x)-a_{C}(t)|
\,dx dt \leq \theta,
\end{equation}
where
$$
a_{C}(t)=\dashint_{C}a (t,x)\,dxds 
\quad (\text{note}\,\,t\,\,\text{and}\,\,ds).
$$

(ii) The  vector-valued $b =(b^{i} )$
admits a representation
$b = b_{M} +b_{B}  $ (``Morrey'' part plus ``bounded'' part) with Borel summands  such that there exist  $  p_{b}\in (d/2, d] $    and
a constant $\hat b_{ M}<\infty$ for which
\begin{equation}
                             \label{6.12.3}
\Big(\dashint_{B } |b_{M}(t,x)|^{p_{b}}\,dx\Big)^{1/p_{b}}\leq \hat b_{M}\rho^{-1}  ,
\end{equation}
whenever $t\in\bR$, $B\in\bB_{\rho}$, and $\rho\leq  r_{b}$,
and there exists a constant $  \|b_{B}\|\in(0,\infty)$
such that
$$
\int_{\bR} \tilde b^{2}_{B} (t)   \,dt\leq    \|b_{B}\|^{2},\quad \tilde b_{B}(t):=
\esssup_{x\in\bR^{d}}|b_{B} (t, x)|.
$$
\end{assumption}
 
Observe that $a_{C}(t)$ is the average of $a(t,x)$
over a ball of radius $\rho$ and, if $a$ is independent of $x$,
the left-hand side of \eqref{6.3.01} is zero. 

Recall that in case $a $ is independent of $t$, we write $a\in VMO$ if 
$a^{\# }_{\rho}\to0$ as $\rho\downarrow 0$. So, our $a$ is ``almost'' in VMO.
It is also worth mentioning that $a\in VMO$ if, for instance,
($a$ is bounded and) $Da\in L_{d}(\mathbb{R}^{d})$
 (see, for instance,
Theorem 10.2.5 in \cite{Kr_08}). An example of such
(uniformly nondegenerate bounded  
highly discontinuous) $a\in VMO$
is given by $2(\delta^{ij})+I_{x \ne0}\zeta(x )\sin(\ln|\ln |x |)$,
where $\zeta$ is any smooth symmetric $d\times d$-matrix valued 
function vanishing for $|x|>1/2$ and
satisfying $|\zeta|\leq 1$.

Introduce
$$
\beta_{b}(t)=\sup_{s}\int_{s}^{s+t}
\tilde b^{2}_{B} (u)   \,du, \quad t\geq0.
$$

\begin{remark}
                   \label{remark 1.5.1}

In the literature  a  very popular condition
on $ b$   is that
$b\in  L_{p ,q  }((0,T)\times\bR^{d})$, that is
\begin{equation}
                            \label{1.5.01}
\|b\|_{ L_{p ,q }((0,T)\times\bR^{d})}=
\Big(\int_{0}^{T}\Big(\int_{\bR^{d}}|b (t,x)|^{p}
\,dx\Big)^{q /p }\,dt\Big)^{1/q }<\infty
\end{equation}
with $p ,q \in[2,\infty]$ satisfying
\begin{equation}
                               \label{9.25.1}
\frac{d}{p }+\frac{2}{q }=1
\end{equation}
see, for instance, \cite{BFGM_19},
\cite{RZ_20}, \cite{RZ_21}, and the references therein.

Observe that, if $p >d$, we can take an
arbitrary constant $\hat N$
and introduce
$$
\lambda(t)=\hat N\Big(\int_{\bR^{d}}
|b (t,x)|^{p }\,dx\Big)^{1/(p -d )},
$$
then for  
$$
b_{M} (t, x)=b  (t, x)I_{|b  (t,x)|\geq \lambda(t)}
$$
and $B\in\bB_{\rho}$  we have
$$
\dashint_{B }|b_{M} (t, x)|^{d}\,dx
\leq \lambda^{d-p }(t)
\dashint_{B }|b  (t, x)|^{p }\,dx
\leq N(d)\hat N^{d-p }\rho^{-d}.
$$
 
Here $N(d)\hat N^{d-p }$ can be made
arbitrarily small if we choose $\hat N$
large enough. In addition, for 
$b_{B} =b -b_{M} $ we have
 $|b_{B} |\leq \lambda $ and
$$
\int_{0}^{T}\lambda^{2}(t)\,dt
=\hat N^{2}\int_{0}^{T}\Big(\int_{\bR^{d}}|b (t,x)|^{p }
\,dx\Big)^{q /p }\,dt<\infty.
$$
This shows that Assumption 
\ref{assumption 6.3.1} (ii) is weaker
than  \eqref{1.5.01}, which is supposed to hold as one of alternative
assumptions
in \cite{RZ_20} and \cite{BFGM_19} if $p>d$.

In case $p=\infty$ and $q=2$,
Assumption 
\ref{assumption 6.3.1} (ii) is the same as in   \cite{BFGM_19},
just take $b_{M}=0$  (this case is not 
considered in   \cite{RZ_20}, \cite{RZ_21}), but, if $p =d$ (and $q =\infty$) our condition is,  basically, weaker than in \cite{RZ_20} (this case is excluded in \cite{RZ_21})
and \cite{BFGM_19}
since we can take $b_{M}=bI_{|b|\geq\lambda}$,
where $\lambda$ is a large constant
and
$$
\int_{\bR^{d} }|b_{M} (t, x)|^{d}\,dx
$$
will be uniformly small if $b(t,\cdot)$
is a continuous $L_{d} $-valued function (one of alternative conditions
in \cite{RZ_20} and \cite{BFGM_19})
or the $L_{d} $-norm of
$b_{M}(t,\cdot)$ is uniformly small  as
in  \cite{BFGM_19}. Our condition on $b$
is satisfied with $p_{b}=d$, for instance,   if
$$
\lim_{\lambda\to\infty}\sup_{[0,T]}\int_{\bR^{d}}|b(t,x)|^{d}I_{|b(t,x)|\geq
\lambda \xi(t) }\,dx <\hat b_{M}^{d} 
$$
for a  function $\xi(t)$
of class $L_{2}([0,T])$.

If $p=d$ an alternative condition
on $b$  in \cite{RZ_20} is that
\begin{equation}
                       \label{1.31.1}
\lambda^{d}|B\cap\{|b(t,\cdot)|>\lambda\}|
\end{equation}
should be sufficiently small uniformly 
for all $B\in\bB_{1}$, $t$, and $\lambda>0$.
It turns out that in this case 
Assumption 
\ref{assumption 6.3.1} (ii) is satisfied (with $b_{B}=0$), any $p_{b}\in (1,d)$, and $r_{b}=1$. This is shown in the following way, where
$B\in\bB_{r},r\leq 1$, $\alpha^{d}=M$,
and $M$ is the supremum of expressions in \eqref{1.31.1}
$$
r^{p_{b}}\dashint_{B}|b(t,x)|^{p_{b}}\,dx
=N(d)r^{p_{b}-d}\Big(\int_{0}^{\alpha/r}
+\int_{\alpha/r}^{\infty}\Big)\lambda^{p_{b}-1}
|B\cap\{|b(t,\cdot)|>\lambda\}|\,d\lambda
$$
$$
\leq Nr^{p_{b}-d}\int_{0}^{\alpha/r}r^{d}
\lambda^{p_{b}-1}\,d\lambda
+Nr^{p_{b}-d}M\int_{\alpha/r}^{\infty}\lambda^{p_{b}-d-1}\,d\lambda
=NM^{p_{b}/d}.
$$

\end{remark}

The case when one has $<$ in place of
$=$
in \eqref{9.25.1} is usually called
subcritical, whereas \eqref{9.25.1}
is a critical case. It turns out that
Assumption 
\ref{assumption 6.3.1} (ii) can be satisfied with
 $p_{b}<d$   and $b_{M}(t,\cdot)\not\in L_{p_{b}+\varepsilon,\loc}$ no matter how small $\varepsilon>0$ is.
In this sense we are dealing with a ``supercritical'' case.  Also note that $\sigma$ is constant in 
\cite{BFGM_19}, \cite{Ki_23},
\cite{RZ_20} and many other papers.
\begin{example}
                        \label{example 9.25.1}

 Take $p \in[d-1,d)$ and take $r_{n}>0$, $n=1,2,...$, such that
the  sum of $\rho_{n}:=r_{n}^{d-p}$ is $1/2$. Let $e_{1}$ be the first
basis vector  and set $b(x)=|x|^{-1}
I_{|x|<1}$, $x_{0}=1$,
$$
x_{n}=1-  2\sum_{1}^{n}r_{i}^{d-p},\quad 
c_{n}=(1/2)(x_{n}+x_{n-1})
$$
$$
  b_{n}(x)=r_{n}^{-1}b\big(r_{n}^{-1}
(x-c_{n}e_{1})\big),\quad b=\sum _{1}^{\infty}b_{n}.
$$
Since $r_{n}\leq 1$ and $d-p\leq 1$,  the supports of $b_{n}$'s are disjoint and
for $q>0$
$$
\int_{B_{1}}b^{q}\,dx=\sum _{1}^{\infty}\int_{\bR^{d}}b_{n}^{q}\,dx=N(d,p)\sum_{1}^{\infty}r_{n}^{d-q}.
$$
According to this we take the $r_{n}$'s so that
the last sum diverges for any $q>p$.
Then observe that for any $n\geq 1$ and any ball $B$
of radius $\rho$
$$
 \int_{B }  b_{n} ^{p}dx \leq N(d) \rho^{d-p} .
$$
Also, if the intersection of $B$ with $\bigcup B_{r_{n}}(c_{n})$
is nonempty, the intersection
 consists of some $B_{r_{i}}(c_{i})$, $i=i_{0},...,i_{1}$, and $B\cap B_{r_{i_{0}-1}}(c_{i_{0}-1})$ if $i_{0}\ne 0  $ and 
$B\cap B_{r_{i_{1}+1}}(c_{i_{1}+1})$. 
In this situation
$$
 \sum_{i=i_{0} }^{i_{1} }
\rho_{i} \leq 2\rho,
$$
and therefore,
$$
 \int_{B }  b  ^{p}\,dx=N(d)\sum_{i=i_{0}}^{i_{1}}r_{i}^{d-p}+\int_{B }  b_{i_{0}-1}  ^{p}dx
+\int_{B }  b_{i_{1}+1}  ^{p}dx \leq N(d)(\rho
+\rho^{d-p}),
$$
where the last term is less than $N(d)\rho^{d-q}$
for $\rho\leq 1$ and this yields  just a different form of \eqref{6.12.3}.
\end{example}
 
\mysection{Existence of solutions}
                          \label{section 3.4.1}

By classical results, if our coefficients are regular enough
(what is assumed in Theorem \ref{theorem 6.5.1}), given a 
$d $-dimensional
Wiener process $w_{t}$, for any $(t,x)\in\bR^{d+1}$ one can uniquely solve equation \eqref{6.15.2}
and the solutions $(t+s,x_{s})$ form a strong Markov process. Sometimes we write $x_{s}
=x_{s}(t,x)$ to emphasize that $x_{s}$
is a solution of \eqref{6.15.2} with given $(t,x)$.

The following theorem shows that, to satisfy 
\eqref{6.5.5},
it suffices to require $\theta$ and $\hat b_{M}$
to be small enough. 
However this is, so far, only for
  regular coefficients.
We use $m_{b}=m_{b}(d,\delta)\in(0,1]$ that
is introduced in Theorem Theorem 3.3 of \cite{Kr_pta}.  

\begin{theorem}
                         \label{theorem 6.5.1}
Suppose that $\sigma$ and $b$
are smooth in $x$  and bounded in $(t,x)$
along with any derivative of any order with
respect to $x$.  
   
Then there exists $\hat b_{M}=\hat b_{M}
(d,\delta,p_{b}, r_{a})>0$,   and there exist $\theta=\theta(d,\delta,  p_{b})>0$, and $\rho_{b}=
\rho_{b}( r_{b} , m_{b},\beta_{b}(\cdot))>0$  such that
  if   Assumption 
\ref{assumption 6.3.1} is satisfied with the above $\theta$, $\hat b_{M} $,     then      
for any $\rho\leq \rho_{b}$,
$C\in\bC_{\rho}$, and
$(t,x)\in\bR^{d+1}$
\begin{equation}
                                 \label{6.5.5}
E_{t,x}\int_{0}^{\tau_{C}}|b(t+s,x_{s})|\,ds
\leq m_{b}  \rho,
\end{equation}
where   $\tau_{C}$
is the first exit time of $(t+s,x_{s})$
from $C$.
\end{theorem}

We prove this theorem after some preparations
done under {\em its assumption that $\sigma$ and $b$
are smooth in $x$  and bounded in $(t,x)$
along with any derivative of any order with
respect to $x$.}
 Set
$$
\cL_{0} u (t,x)=\partial_{t}u (t,x)+ (1/2)a^{ij} (t,x)D_{ij}u (t,x) .
$$
Here is a corollary of Theorem 6.3
of \cite{DK_18}. For $T\in(0,\infty)$ set $\bR^d_{T}=(0,T)\times\bR^{d}$ and fix some $p,q\in(1,\infty)$.

\begin{theorem}
                       \label{theorem 6.3.1}
There exists $\theta =\theta (d,\delta,p,q)>0$
such that if Assumption \ref{assumption 6.3.1}
(i) is satisfied with  this $\theta $, then for any
$f\in  L_{p,q}(\bR^d_{T})$ and $\lambda\geq0$ there exists a unique
$u\in  W^{1,2}_{p,q}(\bR^d_{T})$ such that
$u(T,x)=0$ and $\cL_{0} u -\lambda u =f $ in $\bR^d_{T}$
(a.e.). Moreover, there exists   constants
$(\lambda_{0},N_{0})=(\lambda_{0},N_{0})(d,\delta,p,q, r_{a})$ such that for $\lambda\geq \lambda_{0}$
\begin{equation}
                                 \label{6.3.2}
\|\lambda u,\sqrt\lambda Du, D^{2}u\|_{ L_{p,q}(\bR^d_{T})}\leq N_{0}
\|f\|_{ L_{p,q}(\bR^d_{T})}.
\end{equation}

\end{theorem}

It is perhaps worth saying a few words about
how this theorem is derived from Theorem 6.3
of \cite{DK_18}, which in fact only contains
an a priori estimate for solutions
of $\cL_{0} u -\lambda u =f $ in the {\em whole\/}
$\bR^{d+1}$. By taking smooth $f$ with
compact support in $\bR^{d}_{T}$ one solves
the model equations $\partial_{t}u+\Delta u-\lambda u=f$ by means of a well-known explicit formula,
showing that the solution is smooth bounded and along with any derivative of any order in $t$ and $x$ decreases exponentially fast as $|t|+|x|\to\infty$. Furthermore, this solution vanishes for $t\geq T$. By the a priori estimate in $\bR^{d+1}$ we have \eqref{6.3.2}
for this solution. Then by approximation, we prove that the theorem is true for our model equation and the method of continuity finishes
the job. 

 \begin{remark}
                    \label{remark 11.30.1}
As it is pointed out in \cite{DK_18}, the a priori
estimate we were talking about above is also
true if
we interchange the order of integrations with respect to
$t$ and $x$ in the definition of $ W^{1,2}_{p,q}$. 
However, going to the full operator $\cL$ from
$\cL_{0}$ {\em in our present setting\/} requires $ W^{1,2}_{p,q}$
and $ L_{p,q}$ as we have defined because of \eqref{12.5.4} below.
\end{remark}

Set
$$
\cL_{M}u (t,x)=\partial_{t}u (t,x)+  (1/2) a^{ij} (t,x)D_{ij}u (t,x)+b^{ i}_{M} (t,x)D_{i}u (t,x) ,
$$
\begin{equation}
                                         \label{12.10.4}
\cL u (t,x)=\partial_{t}u (t,x)+ (1/2) a^{ij} (t,x)D_{ij}u (t,x)+b^{ i}  (t,x)D_{i}u (t,x).
 \end{equation}

\begin{theorem}
               \label{theorem 6.3.2}
Suppose that Assumption \ref{assumption 6.3.1} (i) is satisfied with  $\theta$ from Theorem \ref{theorem 6.3.1}.
If $p\in(1,  p_{b}),q\in(1,\infty)$ and
\begin{equation}
                           \label{6.3.5}
N_{0}N_{1}(d,p,p_{b})\hat b_{M} \leq 1/2,
\end{equation}
where $N_{0}=N_{0}(d,\delta,p,q,r_{a})$ is taken from Theorem \ref{theorem 6.3.1} and $N_{1}(d,p,p_{b})$ is specified in the proof,
then  for any
$f\in  L_{p,q}(\bR^d_{T})$ there exists a unique
$u\in  W^{1,2}_{p,q}(\bR^d_{T})$ such that
$u(T,x)=0$ and $\cL_{M}u =f $ in $\bR^d_{T}$
(a.e.). Moreover,  
\begin{equation}
                                 \label{6.3.3}
\|u\|_{ W^{1,2}_{p,q}(\bR^d_{T})}\leq  N 
\|f\|_{ L_{p,q}(\bR^d_{T})},
\end{equation}
where $N=N(d,\delta,p,p_{b},q, r_{a}, r_{b},T)$.

\end{theorem}

Proof. In light of the method of continuity,
it suffices to prove \eqref{6.3.3} as an a priori
estimate. From \eqref{6.3.2} we have
$$
\lambda\|u\|_{ L_{p,q} }+\sqrt{\lambda}\|Du\|_{ L_{p,q} }+\|D^{2}u\|_{ L_{p,q} }
$$
\begin{equation}
                                 \label{6.4.1}
\leq  N_{0} 
\|L_{M}u-\lambda u\|_{ L_{p,q} }+N_{0}\|b^{ i}_{M}  D_{i}v\|_{ L_{p,q} }.
\end{equation}
Here is a crucial estimate which is a corollary
of a Chiarenza-Frasca result.
Since $p<p_{b}\leq d$ we infer from Lemma 3.5
of \cite{1} that with $N_{1}=N_{1}(d,p,p_{b})$ and any $t$
\begin{equation}
                                               \label{12.5.4}
\|\,|b_{M} |\,|Dv |\|_{L_{p}(\bR^{d})}
\leq N_{1}\hat b_{M}\| D^{2}v \|_{L_{p}(\bR^{d})}
+N_{1} r_{b}^{-1}\hat b_{M} \| D v \|_{L_{p}(\bR^{d})}.
\end{equation}
This shows how to estimate the last term in \eqref{6.4.1} and choose $\hat b_{M} $ in order to absorb the term into the left-hand side of \eqref{6.4.1} on the account of also choosing $\lambda$ large enough. Then the usual change
of the unknown function $u$ to $ue^{\lambda t}$
 (the place where the dependence on $T$ appears) yields \eqref{6.3.3}.  The theorem is proved.
\qed

\begin{theorem}
                       \label{theorem 6.4.2}
Assume  
\begin{equation}
                                 \label{6.4.20} 
p\in(1, p_{b}),\quad q\in(1,\infty),\quad\frac{d}{p}+\frac{2}{q}<2,
\end{equation}
(recall that $p_{b}>d/2$)
and suppose that Assumption \ref{assumption 6.3.1} 
(i) is satisfied with  $\theta$ from Theorem \ref{theorem 6.3.1} 
and \eqref{6.3.5} holds. Then
for any $(t,x)\in\bR^{d+1}$, $m=1,2,...$, and Borel nonnegative $f$ on $\bR^{d+1}$ we have
\begin{equation}
                                 \label{6.4.3}
I_{ m}:=E_{t,x}\Big(\int_{0}^{T}f(t+s,x_{s})\,ds\Big)^{m}
\leq m!N^{m}\|f\|^{m}_{ L_{p,q} }
\exp \int_{t}^{t+T} mN(\delta)  \tilde b^{2}_{B}(s)\,ds,
\end{equation}
where   $N$ depends
only on $d,\delta,p,  p_{b},q, r_{a}, r_{b},T$.
\end{theorem}

Proof.   The usual way to treat the additive functionals
of Markov processes  
shows that we may concentrate on $m=1$. Then first assume that $b_{B}\equiv0$ and observe that it suffices to prove 
\eqref{6.4.3} for smooth bounded $f$
with compact support and $(t,x)=0$. 
By classical results for such
an $f$ there exists a regular (of class
$ W^{1,2}_{r,r}(\bR^{d}_{T})$ for any $r$) solution of the problem $\cL u+f=0$ in 
$\bR^{d}_{T}$, $u(T,\cdot)=0$. By It\^o's formula
(applicable to $ W^{1,2}_{d+1,d+1}(\bR^{d}_{T})$-functions) one sees that $I=u(0,0)$ and 
\eqref{6.4.3} follows from \eqref{6.3.3}
and embedding theorems applicable in light of
\eqref{6.4.20}. 

Now still for $b_{B}\equiv0$ we have that
there is a constant $N$ of the same type as in the statement of the theorem such that
for any $m=1,2,...$ we have
\begin{equation}
                                 \label{6.4.4}
E_{t,x}\Big(\int_{0}^{T}f(t+s,x_{s})\,ds\Big)^{m}
\leq m!N^{m}\|f\|^{m}_{ L_{p,q} }.
\end{equation}

Next, we consider the general case.
Denote by $y_{s}$ the solution of 
\eqref{6.15.2} with $b_{M}$ in place of $b$.
Then by Girsanov's theorem
\begin{equation}
                                 \label{6.4.50}
I_{ 1}=E_{t,x} e^{\phi_{T}} \int_{0}^{T}f(t+s,y_{s})\,ds,
\end{equation}
where  
$$
\phi_{ s}= \int_{0}^{s}\gamma_{ u}\,dw_{u}
-(1/2)\int_{0}^{s}|\gamma_{ u}|^{2}\,du ,\quad
\gamma_{ u} = \big (\sigma ^{-1} b_{B} \big)(u,y_{u}). 
$$ 
We are given that $|b_{B} (t+s,y_{s})|\leq
\tilde b _{B}(t+s)$ and therefore, simple and well-known manipulations show that
$$
E_{t,x} e^{2\phi_{T}} \leq E_{t,x}\exp
\int_{0}^{T} N(\delta) |b_{B} (t+s,y_{s})|^{2}\,ds\leq
\exp \int_{t}^{t+T} N(\delta) \tilde b^{2}_{B}(s)\,ds .
$$
This and \eqref{6.4.4} allow  us to obtain
\eqref{6.4.3} with $m=1$ by applying H\"older's
inequality to \eqref{6.4.50}. The theorem is proved.
\qed

\begin{remark}
                              \label{remark 9.7.1}
For $T\in(0, r_{b}^{2}]$ the constant $N$ in \eqref{6.4.3}
can be taken in the form $N_{2}T^{1-(1/2)(d/p+2/q)}$,
where $N_{2}=N_{2}(d,\delta,p,  p_{b},q, r_{a})$ (independent of $ r_{b}$).

To prove this  we may assume that $(t,x)=0 $
and then set 
$$
 \sigma'(t,x)=\sigma(Tt,\sqrt T x),
\quad  b'(t,x)=\sqrt T b(Tt,\sqrt T x).
$$
Since $T\leq1$, Assumption \ref{assumption 6.3.1} is satisfied for $\sigma',b'$ with the same $ r_{a}$
and maybe even smaller $\theta$  and
$\hat b'_{M}$ is the same
as $\hat b_{M}$ and with the corresponding
$ r_{b}'=\bar\rho/\sqrt T\geq 1$.
  Note that condition \eqref{6.3.5}
is satisfied by $\hat b'_{M}$. Then owing to Theorem \ref{theorem 6.4.2} applied in the case of $ r_{b}=1$, after defining
$x'_{s}$ as a unique solution of
$$
x'_{s}=\int_{0}^{s}\sigma'(u,x'_{u})\,dw_{u}
+\int_{0}^{s}b'(u,x'_{u})\,du,
$$
we obtain
\begin{equation}
                             \label{7.11.1}
E_{0,0}\int_{0}^{1}f(T s,\sqrt T x'_{s})\,ds
\leq N  \|f(T\cdot,\sqrt T\cdot)\|_{ L_{p,q} }
\exp \int_{0}^{1}  N(\delta) T \tilde b^{2}_{B}(Ts)\,ds
$$
$$
=NT^{-(1/2)(d/p+2/q)} \|f \|_{ L_{p,q} }
\exp \int_{0}^{T}  N(\delta)   \tilde b^{2}_{B}( s)\,ds,
\end{equation}
where   $N$ depend 
only on $d,\delta,p,  p_{b},q, r_{a}$.

  Observe that the   left-hand side of \eqref{7.11.1}
is  
$$
T^{-1}E_{0,0}\int_{0}^{T}f (  s,\sqrt T x'_{s/T})\,ds,
$$
where, as is easy to see, 
$y_{s}:=\sqrt T x'_{s/T}$ satisfies
$$
y_{s}=\int_{0}^{s}\sigma (u,y_{u})\,dw'_{u}
+\int_{0}^{s}b(u,y_{u})\,du,
$$
with $w'_{s}=\sqrt T w_{s/T}$ which is a Wiener process. Since the coefficients $\sigma$ and $b$ are regular the distribution of $y_{\cdot}$
coincides with that of $x_{\cdot}$,
solution of \eqref{6.15.2} for $t=0,x=0$. Hence,
$$
 T^{-1}E_{0,0}\int_{0}^{T}f (  s,  x _{s })\,ds  \leq NT^{-(1/2)(d/p+2/q)} \|f \|_{ L_{p,q} }
\exp \int_{0}^{T}  N(\delta)   \tilde b^{2}_{B}( s)\,ds
$$ 
and this is what is claimed.
\end{remark}

Observe that for $\rho\leq r_{b}$, $C\in \bC_{\rho}$ and 
$\tau_{C}$ defined as the first exit time of
$(t+s,x_{s})$ from $C$, obviously,
$\tau_{C}\leq \rho^{2}$.
Therefore, Remark \ref{remark 9.7.1}
implies the following.

\begin{corollary}
                   \label{corollary 6.5.1}
For $\rho\leq r_{b}$, $C\in \bC_{\rho}$,  and Borel $f\geq0$
\begin{equation}
                                 \label{6.4.5}
E_{t,x}\int_{0}^{\tau_{C}}f(t+s,x_{s})\,ds
\leq N  \rho ^{2}\dashnorm f\|_{L_{p,q}(C)}  \exp \int_{t}^{t+\rho^{2}}  N(\delta)  \tilde b^{2}_{B}(s)\,ds,
\end{equation}  
where
 $N$ depends
only on $d,\delta,p,  p_{b},q,r_{a}  $.
\end{corollary}
 
 It is convenient to introduce
{\em arbitrary\/} functions 
\begin{equation}
                             \label{7.23.3}
p(d,t),\quad q(d,t),\quad t\in(1,d)
\end{equation}
such that \eqref{6.4.20} holds with $p=p(d,t),
q=q(d,t)$ and any $t$. This is possible
since $p_{b}>d/2$. 

 {\bf Proof of Theorem \protect\ref{theorem 6.5.1}}. Recall that 
$\sigma $ and $b$
are assumed to be smooth  and bounded.  
  Set   $p=p(d,p_{b}) ,q=q(d,p_{b}) $
and  
suppose that Assumption \ref{assumption 6.3.1} (i) is satisfied with   $\theta=\theta(d,\delta,p_{b}):=\theta(d,\delta,p,q)$, where the last $\theta$ is taken from Theorem \ref{theorem 6.3.1}.
Then, as the first restriction, require $\hat b_{M }$ to satisfy \eqref{6.3.5}, which now is a restriction
in terms of only $d,\delta,p_{b}, r_{a} $.
In this situation by \eqref{6.4.5}
  for any $\rho\leq  r_{b}$, $C\in\bC_{\rho}$, and $(t,x)
\in \bR^{d+1}$
$$
E_{t,x}\int_{0}^{\tau_{C}}|b_{M}(t+s, x_{s})|\,ds
\leq   N_{2}\rho^{2}\dashnorm b_{M}\|_{ L_{p,q}(C)}
e^{N(\delta)\beta_{b}(\rho^{2})},
$$
where $N_{2}$ now depends only on $d,\delta,   p_{b},  r_{a} $.
H\"older's inequality shows that
$$
\dashnorm b_{M}\|_{ L_{p,q}(C)}\leq \hat b_{M}\rho^{-1}. 
$$

We subject $\rho_{b}$ to     
\begin{equation}
                                  \label{12.23.01}
  \rho_{b}\leq r_{b},
 \quad\beta_{b}(\rho_{b}^{2}) \leq m_{b}/2\quad(<1).
\end{equation}
Then the second and the final restriction on
$\hat b_{M}$ is
$$
N_{2}\hat b_{M}e^{N(\delta)}\leq m_{b}/2.
$$
In that case for any $\rho\leq  \rho_{b}$, $C\in\bC_{\rho}$, and $(t,x)
\in C$
$$
E_{t,x}\int_{0}^{\tau_{C}}|b_{M}(t+s, x_{s})|\,ds
\leq \rho m_{b}/2,
$$
and, obviously,
$$
\int_{0}^{\tau_{C}}|b_{B}(t+s, x_{s})|\,ds
\leq \int_{0}^{\rho^{2}}\tilde b_{B}(t+s ) \,ds
\leq \rho\beta_{b}(\rho^{2})\leq \rho m_{b}/2.
$$
The combination of these inequalities leads to \eqref{6.5.5}  and
proves the theorem. \qed 

\begin{remark}
                        \label{remark 12.23.1}
Conditions \eqref{12.23.01} show that,
if \eqref{6.12.3} holds with any $\rho>0$
and $b_{B}\equiv0$,
then \eqref{6.5.5} is satisfied with any
$r_{b}\in(0,\infty)$

\end{remark}

In what follows below in this   section 
we suppose that 

{\em Assumption 
\ref{assumption 6.3.1} is satisfied with   $\theta$,   $\hat b_{M}$  from Theorem  \ref{theorem 6.5.1} and no
additional  regularity assumptions
are imposed on $a,b$.}

\begin{theorem}
                              \label{theorem 12.4.1}

(i) There is a probability space 
$(\Omega ,\cF ,P )$,
a filtration of $\sigma$-fields $\cF _{s}\subset \cF $, $s\geq0$,
a process $w _{s}$, $s\geq0$, which is a 
 $d$-dimensional Wiener process
relative to $\{\cF _{s}\}$, and an $\cF _{s}$-adapted
process $x_{s}$ such that 
 (a.s.) for all   $s\geq0$ equation \eqref{6.15.2} holds with $(t,x)=(0,0)$;

(ii) Condition \eqref{6.5.5} is satisfied 
with    $ \rho_{b} $ introduced in Theorem  \ref{theorem 6.5.1} (and $(t,x)=(0,0)$);

 For
$p=p(d,p_{b}) ,q=q(d,p_{b}) $

(iii) for any $T\in(0,\infty)$ and
 any   $m=1,2,...$, there exists a constant $\hat N$
such that for any  Borel nonnegative $f$ on
 $\bR^{d+1}$ we have
\begin{equation}
                                 \label{12.6.1}
E \Big(\int_{0}^{T}f(s,x_{s})\,ds\Big)^{m}
\leq  \hat N\|f\|^{m}_{ L_{p,q} }  ,
\end{equation}

 (iii') in addition, the constant $\hat N$
in \eqref{12.6.1} is expressed as
$$
\hat N=m!N^{m}
\exp \int_{0}^{T} mN(\delta)  \tilde b^{2}_{B}(s)\,ds,
$$
  where   $N$ depends
only on $d,\delta,p,  p_{b},q, r_{a}, r_{b},T$, and if   $T\leq  r_{b}^{2}$, the constant $N$  
can be taken in the form $N_{1}T^{1-(1/2)(d/p+2/q)}$,
where $N_{1}=N_{1}(d,\delta,p,  p_{b},q, r_{a})$ (independent of $ r_{b}$);

  For $p=p(d,p_{b}) ,q=q(d,p_{b}) $ 

(iv) for any  $T\in(0,\infty)$,
and $v\in  W^{1,2}_{p,q}(\bR^{d}_{T}) $,   with probability one for all $t  \in[0,T]$
\begin{equation}
                                          \label{12.5.1}
v(t,x_{t})=v(0,0)+\int_{0}^{t }
\cL v(s,x_{s})\,ds+\int_{0}^{t }D_{i}v(s,x_{s})
\sigma^{ij}(s,x_{s}) \,dw^{j}_{s},
\end{equation}
where the last term is a square-integrable martingale.

\end{theorem}

Proof. Take a nonnegative $\zeta\in C^{\infty}_{0}(\bR^{d})$
with unit integral and for $\varepsilon>0$ set
$\zeta_{\varepsilon}(x)=\varepsilon^{-d}\zeta(x/\varepsilon)$ and use the notation 
$u^{(\varepsilon)}=
u*\zeta_{\varepsilon}$,
where the convolution is performed only with respect to $x$.
It is not hard to check that for each $\varepsilon>0$,
regular
$ \sigma^{(\varepsilon)},
   b^{(\varepsilon)}  $ satisfy Assumption 
\ref{assumption 6.3.1}  with   $\theta$ and $\hat b_{M}$ from Theorem  \ref{theorem 6.5.1}. Therefore, 
in light of Theorem  \ref{theorem 6.5.1}, for
the corresponding Markov processes $(t+s,x^{\varepsilon}_{s})$
all results of \cite{Kr_pta}
are applicable.
In particular, by Corollary 3.10 of \cite{Kr_pta}
for any $\varepsilon,n>0$ and
  $r>s\geq 0$  
\begin{equation}
                                    \label{12.5.2}
E_{0,0}\sup_{u\in[s,r]}|x^{\varepsilon}_{u}-x^{\varepsilon}_{s}|^{ n}
\leq N(  |r-s| ^{ n/2}+|r-s| ^{ n}),
\end{equation}
where $N=N(n, r_{b},d,\delta)$. This implies that the  $P_{0,0}$-distributions of $x^{\varepsilon}_{\cdot}$
are precompact on $C([0,\infty),\bR^{d})$
and a subsequence as $\varepsilon=\varepsilon_{n}
\downarrow 0$ of them converges
to the distribution of a process $x _{\cdot}$ defined
on a probability space (the coordinate process on $\Omega=C([0,\infty),\bR^{d})$ 
with cylindrical $\sigma$-field $\cF$ 
completed with respect to
 $P$, which is
 the limiting distribution of $x^{\varepsilon}_{\cdot}$). 
Furthermore, by Theorem \ref{theorem 6.4.2}
for   for $p=p(d,p_{b}) ,q=q(d,p_{b}) $
(note that $\theta$ was chosen after
$ p(d,p_{b}) , q(d,p_{b}) $ were defined),
and $T\in(0,\infty)$
\begin{equation}
                            \label{3.18.1}
E_{0,0}\Big(\int_{0}^{T}f(t,x^{\varepsilon}_{t})\,dt\Big)^{m}
  \leq m!N^{m}\|f\|^{m}_{ L_{p,q} }
\exp \int_{0}^{T} mN(\delta)  \tilde b^{2}_{B}(s)\,ds
\end{equation}
where $N$ depends
only on $d,\delta,p,  p_{b},q, r_{a}, r_{b},T$.
As usual, this estimate is extended  to $x _{\cdot}$
and yields \eqref{12.6.1}. The form of $N$
in case $T\leq  \rho^{2}_{b}$ is indicated in Remark \ref{remark 9.7.1}. Assertion (ii)
is derived from \eqref{12.6.1} as in Corollary
\ref{corollary 6.5.1}.

Observe that
 estimate \eqref{3.18.1} also implies
that for any bounded Borel $f$ with compact support
\begin{equation}
                                    \label{3.15.5} 
\lim_{\varepsilon\downarrow 0}E _{0,0} \int_{0}^{T}  
f(s,x^{\varepsilon}_{s} )\,ds =
E _{0,0} \int_{0}^{T}  
f(s,x^{0}_{s} )\,ds .
\end{equation}

Now  we prove that assertions (i)  holds for
$x_{\cdot}$.
Estimate \eqref{12.5.2} implies that for any finite $T$
$$
\lim_{c\to\infty}P(\sup_{s\leq T}|x^{0}_{s}|>c)=0,
$$
and estimate   \eqref{12.6.1} 
shows that for any finite $c$
$$
E   \int_{0}^{T}  I_{|x^{0}_{s}|\leq c}
|b (s,x^{0}_{s} )|\,dt <\infty.
$$
Hence, with probability one
$$
\int_{0}^{T}  
|b(s,x^{0}_{s} )|\,dt<\infty.
$$

Next, for $0\leq t_{1}\leq...\leq t_{n}\leq t\leq s$ bounded continuos $\phi(x(1),...,x(n))$
and smooth bounded $u(t,x)$ with compact support by It\^o's formula we have 
$$
E_{0,0}\phi(x^{\varepsilon}_{t_{1}},...,x^{\varepsilon}_{t_{n}})
\Big[u(s,x^{\varepsilon}_{s})-u(t,x^{\varepsilon}_{t})-
\int_{t}^{s}\cL^{\varepsilon}u(r,x^{\varepsilon}_{r})\,dr\Big]=0,
$$
where
$$
\cL^{\varepsilon}u=\partial_{t}u+(1/2)a^{\varepsilon ij}D_{ij}u+
b^{\varepsilon i}D_{i}u,\quad a^{\varepsilon}=(\sigma^{(\varepsilon)})^{2}.
$$

Using \eqref{3.18.1}  and the fact that $u$ has compact support show   that
$$
E_{0,0}\int_{t}^{s}\big|\big(b^{\varepsilon}-b\big)(r,x^{\varepsilon}_{r})|Du(r,x^{\varepsilon}_{r})|
\,dr \to0
$$
as $\varepsilon\downarrow 0$. After that we easily conclude that
$$
E \phi(x^{0}_{t_{1}},...,x^{0}_{t_{n}})
\Big[u(s,x^{0}_{s})-u(t,x^{0}_{t})-
\int_{t}^{s}\cL u(r,x^{0}_{r})\,dr\Big]=0.
$$
It follows that the process
$$
u(s,x^{0}_{s}) -
\int_{0}^{s}\cL u(r,x^{0}_{r})\,dr
$$
is a martingale with respect to the completion
 of $\sigma\{x^{0}_{t}: t\leq s\}$. Referring to a well-known result
  from Stochastic Analysis  
  (see, for instance, Section 2 in \cite{Kr_21_1}) 
   proves assertion (i).

To prove (iv) set
$$
\gamma=\frac{d}{p }+\frac{2}{q },\quad 
\kappa=\frac{\gamma}{\gamma-1}.
$$
Observe that $2>\gamma>1$ ($p <p_{b}\leq d$)
and $\kappa>2$. It follows from Corollary 5.3
of \cite{Kr_22} that
\begin{equation}
                             \label{7.6.3}
\|Dv\|_{L_{\kappa p ,\kappa q }(\bR^{d}_{T})}
\leq N\| v\|_{ W^{1,2}_{ p , q }(\bR^{d}_{T})},
\end{equation}
where (and below) $N$ is independent of $v$. In turn,
\eqref{7.6.3} and \eqref{12.6.1} imply
that
\begin{equation}
                             \label{7.6.4}
E\int_{0}^{T}|Dv(s,x_{s})|^{\kappa}\,ds
\leq N\| v\|^{\kappa}_{ W^{1,2}_{ p , q }(\bR^{d}_{T})}
\end{equation}
and this proves the last statement of 
the theorem.

We prove It\^o's formula \eqref{12.5.1} as usual by taking smooth $v_{n}$ converging to $v$ in
$ W^{1,2}_{ p , q }(\bR^{d}_{T}) $. Since $\gamma<2$,
$v_{n}$ converge uniformly in $ \bR^{d+1}$
in light of the embedding theorems.
In what concerns the integral terms in
\eqref{12.5.1}, estimate \eqref{12.6.1}
shows that the only term of interest is
$$
I_{n}:=E\int_{0}^{T}
|b(s,x_{s}|\,|D(v_{n}-v)(s,x_{s})|\,ds
=I'_{n}+I''_{n},
$$
where, owing to H\"older's inequality and
\eqref{7.6.4},
$$
I'_{n}:=E\int_{0}^{T}
|b_{B}(s,x_{s}|\,|D(v_{n}-v)(s,x_{s})|\,ds
$$
$$
\leq\Big( \int_{0}^{T}\tilde b_{B}^{\gamma}(s)\,ds\Big)^{1/\gamma}\Big(E\int_{0}^{T}|D(v_{n}-v)(s,x_{s})|^{\kappa}\,ds\Big)^{1/\kappa}\to0
$$
as $n\to\infty$, and
$$
I''_{n}:=E\int_{0}^{T}
|b_{M}(s,x_{s}|\,|D(v_{n}-v)(s,x_{s})|\,ds
$$
$$
\leq N\Big(\int_{0}^{T}\|\,|b_{M}(s,\cdot)|\,
|D(v_{n}-v)(s,\cdot)|\,\|_{L_{p }(\bR^{d})}^{q }\,ds\Big)^{1/q }.
$$
By Lemma 3.5 of \cite{1} the integrand in the above term is dominated by a constant times
$$
\hat b_{M}^{q }\|(v_{n}-v)(s,\cdot)\|_{W^{2}_{p }(\bR^{d})}^{q }
$$
and this shows that $I''_{n}\to 0$
  as $n\to\infty$ as well. The theorem is proved. \qed

\begin{corollary}
                    \label{corollary 2.18.2}
Let $x_{\cdot}$ be a process from
Theorem \ref{theorem 12.4.1} and let 
 $p=p(d,p_{b})$, $q=q(d,p_{b}) $. Then for almost any $t>0$ the random vector $x_{t}$
has a density $p_{t}(x)$ such that
$p_{t}\in L_{p'}(\bR^{d})$, $p'=p/(p-1)$, and, for any
  $T\in(0,\infty)$,  
$$
\int_{0}^{T}\|p_{t}\|_{L_{p'}}^{q'}\,dt<\infty,\quad   q'=q/(q-1).
$$
\end{corollary}

Indeed, \eqref{12.6.1} implies that
there is $G(t,x)$ such that 
$$
\int_{0}^{T}\Big(\int_{\bR^{d}}G^{p'}(t,x)\,
dx\Big)^{q'}\,dt<\infty
$$
and
for nonnegative $f$
$$
E  \int_{0}^{T}f(s,x_{s})\,ds=\int_{0}^{T}
\int_{\bR^{d}}f(t,x)G(t,x)\,dxdt.
$$
Obviously, $G(t,x)$ is a good candidate for
$p_{t}(x)$.

\mysection{Weak uniqueness and a Markov
process}
                          \label{section 3.4.2}

Suppose that Assumption   
\ref{assumption 6.3.1} is satisfied with   $\theta$,   $\hat b_{M}$,   from Theorem  \ref{theorem 6.5.1} (and no
additional  regularity assumptions
are imposed on $a,b$).
\begin{lemma}
                          \label{lemma 12.7.1}
Take the process $x_{s}$ from Theorem \ref{theorem 12.4.1}
and for $n\geq0$ set
$$
\gamma_{n,s} =I_{(n,\infty)}(|b_{B}(s,x_{s})|)\big(\sigma ^{-1}b_{B} \big)(s,x_{s}),
$$
$$
\psi_{n,s}=-\int_{0}^{s}\gamma_{n,u}\,dw_{u}
-(1/2)\int_{0}^{s}|\gamma_{n,u}|^{2}\,du.
$$
Then for any $T \in(0,\infty)$
$$
E\sup_{s\leq T}\big|\exp
\psi_{n,s}-1  \big|\leq   N\big(\int_{0}^{T} \tilde b_{  B,n}^{2}(s)\,ds\Big)^{1/2}\exp\int_{0}^{T}N\tilde b^{2}_{B}(s)\,ds,
$$
where $N=N(\delta)$ and $\tilde b_{ B,n}^{2}(s)=
I_{(n,\infty)}(\tilde b_{B}(s) )\tilde b^{2}_{B}(s)$.
 
\end{lemma}

Proof. We have   
$$
\exp
\psi_{n,s}-1 =-\int_{0}^{s}\exp
\psi_{n,u}\gamma_{n,u}\,dw_{u}.
$$
Hence, by Doob's inequality
$$
\Big(E\sup_{s\leq T}\big|\exp
\psi_{n,s}-1  \big|\Big)^{2}\leq 4E\int_{0}^{T}
|\gamma_{n,s}|^{2}\exp\big(2
\psi_{n,s})\,ds.
$$
Here $|\gamma_{n,s}|\leq N\tilde b_{  B,n}(s)$,
  $\tilde b_{B}(s)$ is a deterministic function
of class $L_{2}(\bR)$ and
$$
E\exp\big(2
\psi_{n,s})=EI^{1/2}_{n,s}  \Big(\exp\int_{0}^{s}6|\gamma_{n,u}|^{2}\,du \Big)^{1/2}
$$
where
$$
I _{n,s}=\exp\Big(-4\int_{0}^{s}\gamma_{n,u}\,dw_{u}
-8\int_{0}^{s}|\gamma_{n,u}|^{2}\,du\Big).
$$
Since $EI_{n,s}\leq1$,
$$
E\int_{0}^{T}
|\gamma_{n,s}|^{2}\exp\big(2
\psi_{n,s})\,ds\leq NE\int_{0}^{T}
\tilde b_{ B,n}^{2}(s)\,ds \exp\int_{0}^{T}N\tilde b^{2}_{B}(s)\,ds.
$$
This   proves
the lemma. \qed

For $n\geq 0$ introduce
$$
\cL^{n}v(t,x)=\partial_{t}v (t,x)+ a^{ij} (t,x)D_{ij}v (t,x)+b^{i}_{M}(t,x)D_{i}v(t,x)
$$
$$
+I_{[0,n]}(|b _{B}(t,x)|)b^{i}_{B}(t,x)D_{i}v(t,x).
$$

\begin{theorem}
                \label{theorem 7.6.4}
Let  $p=p(d,p_{b}) ,q=q(d,p_{b}) $ 
(see \eqref{7.23.3})   and $T\in(0,\infty)$.
Let $c(t,x)$ be a real-valued bounded Borel
function on $\bR^{d+1}$. Then
for any $n\geq0$ and
$f\in  L_{p ,q }(\bR^d_{T})$ there exists a unique
$v\in  W^{1,2}_{p ,q }(\bR^d_{T})$ such that
$v(T,x)=0$ and $\cL^{n}v+cv =f $ in $\bR^d_{T}$
(a.e.).
\end{theorem}

This theorem for $n=0$ and $c\equiv0$
is a particular case of Theorem \ref{theorem 6.3.2}. Adding bounded lower order terms
does not affect much the proof there.

\begin{remark}
                   \label{remark 7.7.1} 
In Theorem \ref{theorem 7.6.4} the terminal condition is zero. However, it also holds
if it is any $v_{T}\in W^{2}_{p }$.
Indeed, for $t\in  [T,2T]$ define $f(t,x)=
(\partial_{t}+\Delta)\big(v_{T}(x)(2-t/T)\big)$
extend $\cL^{n}$ on $[T,2T]$ as $\partial_{t}+\Delta$ and solve with the new $\cL^{n}$
the equation $\cL^{n}v+cI_{t<T}v=f$ on $[0,2T]$
with zero terminal condition. Then by uniqueness
 $v(t,x)=v_{T}(x)(2-t/T)$ on $[T,2T]$
and $v$ solves the original equation on $[0,T]$ with
the terminal data $v_{T}$.
\end{remark}

\begin{theorem}[Conditional weak uniqueness]
                                    \label{theorem 12.6.1}
In all situations when  
 for $p=p(d,p_{b}) ,q=q(d,p_{b}) $ 
the assertions (i) and (iii)  of 
Theorem \ref{theorem 12.4.1},
with $m=1$   and $\hat N$ 
independent of $f$, or   the assertions (i)
and (iv) of Theorem \ref{theorem 12.4.1}, with the last term in \eqref{12.5.1}
being a local martingale, 
 hold true,
the finite dimensional distributions of 
all such solutions of \eqref{6.15.2}
with fixed $(t,x)$
are the same.  

\end{theorem}

Proof. We may assume that
$(t,x)=(0,0)$. As we have seen in the proof of Theorem 
\ref{theorem 12.4.1} assertion (iii) implies (iv).
Therefore we assume that assertions (i)
and (iv) hold true.
  Take a Borel bounded $c$ on $\bR^{d+1}$
with compact support, take a smooth
function $f$ on $\bR^{d}$ with compact support,
take
  $T,n>0$, and let $v_{n}$ be the solution of the equation $\cL^{n}v_{n}+cv_{n}=0$ in $\bR^{d}_{T}$ with   terminal condition $v_{n}(T,\cdot)=f$. Also take $\gamma_{n,s}$ from Lemma
\ref{lemma 12.7.1} 
and introduce
$$
\phi_{n,s}=-\int_{0}^{s}\gamma_{n,u}\,dw_{u}
-(1/2)\int_{0}^{s}|\gamma_{n,u}|^{2}\,du
+\int_{0}^{s}c(u,x_{u})\,du.
$$
By It\^o's formula, applied to $v_{n}(s,x_{s})\exp
\phi_{n,s}$,   we obtain for $s\leq T$
$$
v_{n}(s,x_{s})e^{\phi_{n,s}}=v_{n}(0,0)
+\int_{0}^{s}
e^{\phi_{n,u}}\big(\sigma^{ik}D_{i}v_{n}(u,x_{u})
+v_{n}(u,x_{u})\gamma^{k}_{n,u}\big)\,dw^{k}_{u}.
$$
It follows that for a sequence of stopping times
$\tau_{m}\uparrow T$ (localizing the stochastic integral
above)

\begin{equation}
                               \label{7.6.7}
Ev_{n}(\tau_{m},x_{\tau_{m}})e^{\phi_{n,\tau_{m}}}=v_{n}(0,0)
\end{equation}
Notice that   $v_{n}$ is a bounded continuous function. Hence, by the dominated convergence theorem
and Lemma \ref{lemma 12.7.1} we infer from \eqref{7.6.7} that
$$
Ef(x_{T})e^{\phi_{n,T}}=v_{n}(0,0).
$$
By sending $n\to\infty$ and using Lemma \ref{lemma 12.7.1}
one more time we find that the limit of $v_{n}(0,0)$ exists and
\begin{equation}
                                \label{7.7.4}
Ef(x_{T})\exp\Big(\int_{0}^{T}c(u,x_{u})\,du
\Big)\,ds=\lim_{n\to\infty}v_{n}(0,0).
\end{equation}

Since the last limit is independent of
what solution $x_{s}$ we take and we have the arbitrariness
in $c$ and $f$,  we conclude that all
restricted as in the statement of the theorem
 solutions
have the same finite-dimensional distributions
  (one can find a detailed proof of
this conclusion, for instance, in \cite{KW_11}).
 
The theorem is proved. \qed

 \begin{remark}
                             \label{remark 6.2.1}
 Deriving analytical conditions on $a,b$ guaranteeing
 the weak uniqueness and allowing $b$ to be {\em
 singular\/}
 is a very challenging and hard problem.
 We know of only one case mentioned in the next remark
 when this was achieved.
 \end{remark}

\begin{remark}    
                \label{remark 10.30.10}
In \cite{RZ_20} the conditional weak uniqueness is proved under condition \eqref{1.5.01}
 with $p,q$ satisfying \eqref{9.25.1}. It turns out that
there are situations when Assumption
\ref{assumption 6.3.1} (ii) is satisfied
and
there are no $p,q$ satisfying \eqref{9.25.1} such that \eqref{1.5.01} holds. 
  For instance,
take $b(t,x)$ such that $|b|=cf$, where 
the constant $c>0$ and
$f=(|x|+\sqrt{|t|})^{-1}I_{|x|<1,|t|<1}$.
Then   for $p,q\in[1,\infty]$,
$d/p+2/q=1$, $p>d$, and small $t>0$
\begin{equation}
                         \label{1.28.6}
\int_{B_{1}}\frac{1}{(|x|+\sqrt{|t|})^{p}}\,dx
=|t|^{(d-p)/2}\int_{B_{1/\sqrt t}}\frac{1}{(|x|+1)^{p}}\,dx
\end{equation}
is of order $|t|^{(d-p)/2}$, whose ($q/p$)th power
is $|t|^{-1}$   which is not integrable near zero. In case $p=d$, $q=\infty$ the integral in \eqref{1.28.6} tends to infinity
as $t\downarrow 0$, so \eqref{1.5.01}
 fails again.   The same happens if $p=\infty$.

On the other hand, since $|f|\leq1/|x|$,  for any $p\in[1,d)$, if $|x|\leq 2\rho$
\begin{equation}
                         \label{11.1.1}
\dashnorm f(t,\cdot)\|^{p}_{L_{p}(B_{\rho}(x))}\leq N\rho^{-d}\int_{B_{3\rho} }
\frac{1}{ |y| ^{p }}
\,dy=N\rho^{-p},
\end{equation}
and if $|x|> 2\rho$, then $f(t,\cdot)\leq \rho^{-1}$
on $B_{\rho}(x)$ and the inequality between the extreme terms of \eqref{11.1.1} holds again.
Hence, if $c$ is sufficiently small the  above $b$
satisfies Assumption
\ref{assumption 6.3.1} (ii) and \eqref{6.15.2}
has   solutions which are conditionally unique, say if $\sigma$ is constant.
Recall that our $\sigma$ is not necessarily constant.

Actually,  the fact that it is (unconditionally) weakly unique 
follows from the fact, which we will prove elsewhere,
saying that for weak uniqueness it suffices
to have some $p_{0},q_{0}\in[1,
\infty)$ such that  $p_{0}\geq q_{0}$,  $d/p_{0}+1/q_{0}=1$
and for any $C\in \bC_{\rho}$
 to have
\begin{equation}
                    \label{1.30.1}
\dashnorm f\| _{L_{p_{0},q_{0}}(C) }
\leq N\rho^{-1}
\end{equation}
with $N$ independent of $C$.

That \eqref{1.30.1} holds for {\em any\/} $p_{0},q_{0}$
described above is shown as follows. If $|x|+\sqrt{|t|}\leq 3\rho $ and $p_{0}\geq q_{0}$ ($p_{0}\geq d+1$), then for  $C=C_{\rho}(t,x)
\subset (-25\rho^{2},25\rho^{2})\times B_{5\rho}$ the $q_{0}$-th power of the left-hand side of \eqref{1.30.1} is dominated by
$$
  \rho^{- q_{0}-1}\int_{-25\rho^{2}}^{25\rho^{2}}
\Big(\int_{\bR^{d}}\frac{1}{(|y|+\sqrt{|s|})^{p_{0}}}
\,dy\Big)^{q_{0}/p_{0}}\,ds
=N\rho^{-q_{0}}.
$$
 In case  $|x|+\sqrt{|t|}\geq 3\rho $, we have $f\leq N/\rho$
on $C_{\rho}(t,x) $ and $\dashnorm f\|_{L_{p_{0},q_{0}}(C_{\rho}(t,x))}\leq N/\rho$ again.

  By the way,
G. Zhao (\cite{Zh_20_1}) gave an example showing
that, if in condition \eqref{6.12.3}   we replace $\rho^{-1}$ with $\rho^{-\alpha}$, $\alpha>1$, the weak uniqueness
may fail even in the time homogeneous case
and unit diffusion.
\end{remark}  

By changing the origin we can apply Theorem \ref{theorem 12.4.1}
to prove the solvability of \eqref{6.15.2}
with any initial data $(t,x)$ and get  solutions with the
properties as in Theorem \ref{theorem 12.4.1}.
For such a solution denote by $P_{t,x}$ the distribution of $(\sft_{s},x_{s}),s\geq0$, ($\sft_{s}=t+s$)
on the Borel $\sigma$-field $\cF$ of $\Omega=C([0,\infty),\bR^{d+1})$. For $\omega=(\sft_{\cdot},x_{\cdot})\in \Omega$
set $(\sft_{s},x_{s})(\omega)=(\sft_{s},x_{s})$.
Also set $\frN_{s}=\sigma \{(\sft_{t},x_{t}),t\leq s\}$.

\begin{theorem}
                               \label{theorem 12.6.2}
The process
$$
X=\{(\sft_{\cdot},x_{\cdot}),\infty,\frN_{t},P_{t,x}\}
\quad (\infty\,\, \text{is the lifetime})
$$
is strong Markov with  
strong  Feller resolvent and satisfies \eqref{6.5.5}  
with    $r_{b}$ introduced in Theorem  \ref{theorem 6.5.1}.
\end{theorem}

Proof. Take $v_{n}$ from the proof of Theorem
\ref{theorem 12.6.1} with $c=0$. 
By It\^o's formula for $(t,x)\in \bR^{d}_{T}$
and $ 0\leq r\leq T-t$ we obtain that with
$P_{t,x}$-probability one for all $s\in[r ,T-t]$
\begin{equation}
                                            \label{12.6.4}
v_{n}(\sft_{s} ,x_{s})e^{\phi_{n,s}}=v_{n}(\sft_{r },x_{r })
^{\phi_{n,r}}
+\int_{r }^{s}
e^{\phi_{n,u}}\big(\sigma^{ik}D_{i}v_{n}(\sft_{u},x_{u})
+v_{n}(\sft_{u},x_{u})\gamma^{k}_{n,u}\big)\,dw^{k}_{u},
\end{equation}
where
$$
\gamma_{n,s} =I_{(n,\infty)}(|b_{B}(\sft_{s},x_{s})|)\big(\sigma ^{-1}b_{B} \big)(\sft_{s},x_{s}),
$$
$$
\phi_{n,s}=-\int_{0}^{s}\gamma_{n,u}\,dw_{u}
-(1/2)\int_{0}^{s}|\gamma_{n,u}|^{2}\,du.
$$

From \eqref{12.6.4} with $r=0$ we obtain
$$
E_{t,x}f(x_{T-t})e^{\phi_{n,T-t}}=v_{n}(t,x ),
$$
which along with Lemma \ref{lemma 12.7.1}
show  that $v_{n}$ are uniformly bounded and
\begin{equation}
                                            \label{12.6.6}
E_{t,x}f(x_{T-t})=\lim_{n\to\infty}v_{n}(t,x),
\end{equation}
which, in particular, implies that $E_{t,x}f(x_{T-t})$
is a  Borel function of $(t,x)$.
It is also, obviously, continuous with respect to $T$.
Hence, it is a Borel functions of $(t,x,T)$. This 
fact is obtained for
sufficiently regular $f$, but by usual measure-theoretic
arguments it carries over to all Borel bounded $f$.
It follows that $E_{t,x}f(x_{T})$ is a Borel
function on $\bR^{d+1}$ for any $T>0$ and Borel bounded
$f$. Then $E_{t,x} \big(\phi(t+T)f(x_{T})\big)$ is also Borel
if $\phi$ is Borel and again usual measure-theoretic
arguments allow us to conclude that 
$E_{t,x}  f(\sft_{T},x_{T}))$ is Borel for any
Borel bounded $f$ given on $\bR^{d+1}$.

Then take $0\leq r_{1}\leq...\leq r_{m}=r$
and a  bounded continuous  function $\zeta\big(x(1),...,x(m)\big)$
on $\bR^{md}$ and conclude from \eqref{12.6.4} that
for any stopping time $\tau$ localizing the stochastic
integral in \eqref{12.6.4} 
$$
E_{t,x}I_{r\leq\tau}\zeta(x_{r_{1}\wedge\tau},...,x_{r_{m}\wedge\tau})v_{n}(\sft_{\tau\wedge(T-t)} ,x_{\tau\wedge(T-t)})e^{\phi_{n,\tau\wedge(T-t)}}
$$
$$
=E_{t,x}I_{r\leq\tau}\zeta(x_{r_{1}\wedge\tau},...,x_{r_{m}\wedge\tau})v_{n}(\sft_{r },x_{r })e
^{\phi_{n,r}}.
$$
Sending $\tau\uparrow T-t$, then $n\to\infty$ and using
\eqref{12.6.6} we conclude
$$
E _{t,x}\zeta(x_{r_{1} },...,x_{r_{m} })f(x_{T-t})=
E_{t,x} \zeta(x_{r_{1} },...,x_{r_{m} })E_{\sft_{r},x_{r}}
f(x_{T-(t+r)}).
$$
The arbitrariness of $r_{i}$'s and $\zeta,f$ proves that
$X$ is a Markov process.  

To prove that it is strong Markov it suffices to prove
that its resolvent $R_{\lambda}$ is Feller. We are going to show that, actually, it is strong Feller, that is maps bounded Borel functions into bounded continuous ones.
The resolvent is given as an integral with respect to
$t$ over $(0,\infty)$. Restricting these integrals to
$[0,T]$ we obtain operators which converge to $R_{\lambda}$
strongly in the uniform norm as $T\to\infty$.  Therefore, to prove
that $R_{\lambda}$ is  
strong  Feller it suffices to show that
for any $T\in(0,\infty)$ and bounded continuous $f$
on $\bR^{d+1}$ the function
\begin{equation}
                                           \label{12.7.6}
E_{t,x}\int_{0}^{T}f(\sft_{s},x_{s})\,ds
\end{equation}
is bounded and continuous. The boundedness is obvious.
Furthermore, in light of estimate
\eqref{12.5.2}, the process $x_{s}$
reaches   the areas far from the initial point with small probability.
Therefore, we need to concentrate only on $f$
with compact support.

Fix such an $f$, take  $p=p(d,p_{b})
 ,q=q(d,p_{b}) $, and for $n\geq0$ denote by $v_{n}$
the solution of class $ W^{1,2}_{p,q}(\bR^{d}_{T})$
of the equation $\cL^{n}v_{n}+f=0$ with zero terminal data.
By It\^o's formula for $(t,x)\in \bR^{d}_{T}$
  we obtain that with
$P_{t,x}$-probability one for all $s\in[0 ,T-t]$
$$
v_{n}(\sft_{s} ,x_{s})e^{\phi_{n,s}}=v_{n}(t,x )
-\int_{0 }^{s}
e^{\phi_{n,u}}f(\sft_{u},x_{u})\,du
$$
$$
+\int_{0 }^{s}
e^{\phi_{n,u}}\big(\sigma^{ik}D_{i}v_{n}(\sft_{u},x_{u})
+v_{n}(\sft_{u},x_{u})\gamma^{k}_{n,u}\big)\,dw^{k}_{u}.
$$
 As before, by using localizing stopping times, we get
\begin{equation}
                                            \label{12.7.7}
E_{t,x}\int_{0 }^{T-t}
e^{\phi_{n,s}}f(\sft_{s},x_{s})\,ds=v_{n}(t,x ).
\end{equation}
By virtue of Lemma \ref{lemma 12.7.1}
$$
\Big|E_{t,x}\int_{0 }^{T-t}
e^{\phi_{n,s}}f(\sft_{s},x_{s})\,ds-
E_{t,x}\int_{0 }^{T-t}
 f(\sft_{s},x_{s})\,ds\Big|
$$
$$
\leq 
\sup|f|T
E_{t,x}\sup_{s\leq T}\big|e^{\phi_{n,s}}-1\big| 
\leq N\big(\int_{0}^{T} \tilde b_{ B,n}^{2}(s)\,ds\Big)^{1/2},
$$
where $N$ is independent of $n,t,x$. It follows that the left-hand sides 
of \eqref{12.7.7}, which are continuous due to \eqref{12.7.7},
converge as $n\to\infty$ uniformly to
$$
E_{t,x}\int_{0 }^{T-t}
 f(\sft_{s},x_{s})\,ds,
$$
which is therefore continuous, and this brings the proof of the theorem
to an end. \qed

\mysection{About strong solutions}
                       \label{section 2.28.1}

Here we are dealing with the equation
\begin{equation}
                         \label{2.28.2}
x_{t}=\int_{0}^{t}\sigma(s,x_{s})\,dw_{s}
+\int_{0}^{t}b(s,x_{s})\,ds,
\end{equation}
where $b$ is as before, $\sigma =(\sigma^{ik})$ is Borel with values in the set of $d\times d_{1}$-matrices ($d_{1}\geq d$), and $w_{t}
=(w^{1}_{t},...,w^{d_{1}})$ is a Wiener process
on a probability space.
We suppose that the matrix-valued $a =(a^{ij})=\sigma\sigma^{*}$   takes its values
in $\bS_{\delta}$, that is the set of $d\times d $
symmetric matrices with eigenvalues between $\delta$ and $\delta^{-1}$. Take an $r_{\sigma}>0$.
\begin{theorem}
                         \label{theorem 2.28.1}

Suppose that Assumption   
\ref{assumption 6.3.1} is satisfied with   $\theta$,   $\hat b_{M}$  from Theorem  \ref{theorem 6.5.1}.
Then there could be only one strong solution
of \eqref{2.28.2} on the given probability space for which 
the assertions (i) and (iii)  (without (iii')),
with $m=1$, or (i)
and (iv), with the last term in \eqref{12.5.1}
being a local martingale, 
of Theorem \ref{theorem 12.4.1} hold true.
  More generally, if on a given probability space
there is
at least one strong solution  
of \eqref{2.28.2}, for which 
the assertions (i) and (iii)  
(without (iii')),
with $m=1$, or (i)
and (iv), with the last term in \eqref{12.5.1}
being a local martingale, 
of Theorem \ref{theorem 12.4.1} hold true,
and if there is another solution with the same finite dimensional
distributions as the strong one, then
it coincides with the strong solution.
 
\end{theorem}

Proof. One is tempted to refer
to the result of A. Cherny \cite{Ch_02} saying that weak uniqueness and strong existence imply
the uniqueness of strong solutions. However,
in this result one needs unconditional
weak uniqueness which we do not know how to prove
in the general case. Therefore, we proceed
differently still using the idea from \cite{Ch_02}.

Define $\sigma=(\sigma^{k})=(\sigma^{ik})$ and let   $\tau=\sigma^{*}\sigma$. This is a symmetric  nonnegative definite matrix and the following
is well defined
$$
\Sigma=\lim_{\varepsilon\downarrow 0}
\tau(\tau+\varepsilon I)^{-1},
$$
where $I$ is the $d_{1}\times d_{1}$ identity matrix. As is easy to see by using the diagonal forms,
$\Sigma^{2}=\Sigma$, $\Sigma\tau=\tau$, and ($\text{tr}\,AB=
\text{tr}\,BA$)
$$
\text{tr}\,( \Sigma \sigma^{*}-\sigma^{*})
(\sigma \Sigma-\sigma)=
\text{tr}\,(\Sigma \sigma^{*}\sigma\Sigma
-\sigma^{*}\sigma\Sigma-\Sigma\sigma^{*}\sigma
+\sigma^{*}\sigma)
=\text{tr}\,(-\Sigma\tau+\tau)=0,
$$
so that $\sigma\Sigma=\sigma$.
 Next, by extending our probability
space, if necessary, we suppose that we are also
given a $d_{1}$-dimensional Wiener process
$\bar w_{s}$ independent of $w_{t}$. Define
$$
\xi_{s}=\int_{0}^{s}\Sigma(u,x_{u})\,d\bar w_{u}+\int_{0}^{s}\big(I-\Sigma(u,x_{u})\big)\,
dw_{u}.
$$
An easy application of the L\'evy theorem shows that $\xi_{s}$ is a $d_{1}$-dimensional Wiener process.

The crucial step is to prove that
the processes $x_{\cdot}$ and $\xi_{\cdot}$
are independent on $[0,T]$ because
(dropping arguments $(s,x_{s})$)
$$
dx^{i}_{s}d\xi^{k}_{s}=\sigma^{ir}\,dw^{r}_{s}(
\delta^{kn}-\Sigma^{kn})\,dw^{n}_{s}
=\sigma^{ir} (
\delta^{kr}-\Sigma^{kr})\,ds
=(\sigma^{ik} 
 -\sigma^{ir}\Sigma^{rk})\,ds=0.
$$

To do that, take two bounded Borel
functions $c'$ and $c''$ with compact support on $\bR^{d+1}$ and $\bR^{d_{1}+1}$, respectively, and take a smooth
compactly supported function $f$
on $\bR^{d}$. Then for $n>0$ define
$v'_{n}$ as a solution of $\cL^{n}v'_{n}
+c'v'_{n}=0$ on $[0,T]\times \bR^{d}$ and $v''$ as a solution
of $(\partial_{t}+(1/2)\Delta +c'')v''=0$
$[0,T]\times \bR^{d_{1}}$ with terminal conditions
$v'_{n}(T,\cdot)=f$, $v''(T,\cdot)=1$.

By It\^o's formula applied to 
$$
v'_{n}(s,x_{s})v''(s,\xi_{s})e^{\phi_{n,s}},
$$
where
$$
\phi_{n,s}=-\int_{0}^{s}\gamma_{n,u} \,dw _{u}
-(1/2)\int_{0}^{s}|\gamma_{n,u}|^{2}\,du
+\int_{0}^{s}[c'(u,x_{u})+c''(u,\xi_{u})]\,du,
$$
$$
\gamma_{n,s}=I_{[n,\infty)}(|b_{B}(s,x_{s})|)\big(\sigma^{*}a^{-1}b\big)(s,x_{s}),
$$
we get
$$
f( x_{T}) e^{\phi_{n,T}}=
v' _{n}(0,0)v''(0,0)
$$
$$
+\int_{0}^{T}e^{\phi_{n,u}}\big[v''(u,\xi_{u}) \sigma^{ik}D_{i}v'_{n}(u,x_{u})\,
dw^{k}_{u}+v'_{n}(u,x_{u})D_{\xi^{i}}
v''(u,\xi_{u})\,d\bar w^{i}_{u}\big].
$$
By taking expectations and arguing as in the proof of 
\eqref{12.6.6},
we see that  
\begin{equation}
                   \label{7.7.2}
Ef(x_{T})\exp\Big(\int_{0}^{T}c'(s,x_{s})\,ds\Big)\exp\Big(\int_{0}^{T}c''(s,\xi_{s})\,ds\Big)=\lim_{n\to\infty}v^{'}_{n}(0,0)
v^{''}(0,0).
\end{equation}
The limit here we find from \eqref{12.6.6} and conclude that
$$
Ef(x_{T})\exp\Big(\int_{0}^{T}c'(s,x_{s})\,ds\Big)\exp\Big(\int_{0}^{T}c''(s,\xi_{s})\,ds\Big)
$$
$$
=Ef(x_{T}) \exp\Big(\int_{0}^{T}c'(s,x_{s})\,ds\Big)
v''(0,0).
$$

After that the arbitrariness of $f$
shows that
$$
E \exp\Big(\int_{0}^{T}c'(s,x_{s})\,ds\Big)\exp\Big(\int_{0}^{T}c''(s,\xi_{s})\,ds\Big)
$$
$$
=E \exp\Big(\int_{0}^{T}c'(s,x_{s})\,ds\Big)
v''(0,0).
$$
By taking $c'=0$ we identify $v''(0,0)$
and then the arbitrariness of $c'$ and $c''$ proves that $x_{\cdot}$ and $\xi_{\cdot}$ are independent indeed.

Then we observe that
$$
w_{s}=I_{s} 
+\int_{0}^{s}\big(I-\Sigma(u,x_{u})\big)\,d\xi_{u},
$$
where
$$
I_{s}=\int_{0}^{s}\Sigma(u,x_{u})\,dw_{u}=
\lim_{\varepsilon\downarrow 0}
\int_{0}^{s}(\tau(u,x_{u})+\varepsilon I)^{-1}
\sigma^{*}(u,x_{u})\,dm_{u},
$$
$$
m_{s}=\int_{0}^{s}\sigma (u,x_{u})\,dw_{u}
=x_{s}-\int_{0}^{s}b(u,x_{u})\,du.
$$
We see that $I_{s}$ is a functional of $x_{\cdot}$, so that the distribution
of $I_{\cdot}$ are defined uniquely. Since
the Wiener process
$\xi_{\cdot}$ is independent of $x_{\cdot}$,
the conditional distribution of $w_{\cdot}$
given $x_{\cdot}$ and
the joint distribution of $(w_{\cdot},x_{\cdot})$ are unique.

The remaining part of the proof 
  can now follow exactly as in \cite{Ch_02}. The theorem is proved.\qed

   Concerning the {\em existence\/} of strong solutions
 we have a conjecture before which we introduce
 appropriate assumptions. This conjecture is quite 
 challenging even if $b\equiv0$. 

\begin{assumption}
                     \label{assumption 6.3.10}
(i)  Assumption \ref{assumption 6.3.1} is satisfied with   $\theta$,   $\hat b_{M}$  from Theorem  \ref{theorem 6.5.1},
 $d\geq 3$, $p_{b}\in ( 2   ,d]$  
and (recall \eqref{7.23.3}) $ p(d,p_{b})
\geq 2, p(d,p_{b})>d/2$.

(ii) For any $t$, $\sigma(t,\cdot)\in W^{1}_{1,\loc} $ and the tensor-valued $D\sigma$ 
admits a representation
$D\sigma=D\sigma_{M}+D\sigma_{B}$ 
with Borel summands (``Morrey part'' of
$D\sigma$ plus the ``bounded part'') such that
there exist    a
finite constant $\widehat{D\sigma}_{M}$
and $p_{D\sigma}\in(2,d]$
for which
$$
\Big(\dashint_{B } |D\sigma_{M}(t,x)|^{p_{D\sigma}}\,dx\Big)^{1/p_{D\sigma}}\leq    \rho^{-1}\widehat{D\sigma}_{M}   ,
$$
whenever $t\in\bR$, $B\in\bB_{\rho}$, and $\rho\leq  r_{\sigma}$, 
and there exists a constant $\| D\sigma_{B}\|
 \in(0,\infty) $
such that
$$
\int_{\bR}| \widetilde{D\sigma}_{B} (t)|^{2}   \,dt\leq\|D\sigma_{B}\|^{2},\quad 
\widetilde{D\sigma}_{B} (t) :=
\esssup_{x\in\bR^{d}}|D\sigma_{B} (t, x)|;
$$

\end{assumption}

{\em Conjecture. There exists $\beta=\beta(d,d_{1},\delta,   p_{D\sigma},  p_{b})>0$
such that if 
$$
\widehat{D\sigma}_{M},\hat b_{M}\leq \beta,
$$
 then equation \eqref{2.28.2} has a strong solution
 for which assertions (iii) and (iv) of
 Theorem \ref{theorem 12.4.1} hold true. Thus, it is
 unique in the class of (weak or strong) solutions with the same
 finite-dimensional distributions. 
 
 What follows below is to explain how we came
 to this conjecture and a hypothetical way to prove it.
 By the way, the author had a few failed attempts
 while following that way.}

As in \cite{1} we approximate the coefficients
of \eqref{2.28.2} by smooth functions and
deal with representations of the solutions of approximating
equations. Therefore, we additionally
 suppose that $\sigma $ and $b$
are {\em smooth in $x$\/}  and bounded in $(t,x)$
along with any derivative of any order with
respect to $x$.

Denote by $C^{\infty}_{3} $ the set
of infinitely differentiable functions on $\bR^{d}$   which have bounded partial derivatives
(of any order) summable to the 3rd power
over $\bR^{d}$. Here 3 can be replaced
with any number $>2$ and appears in order
to be able to refer to some results
from \cite{Kr_99}.

For $n=1,2,...$ and $r\in(0,\infty)$
introduce 
$$
\Gamma^{n}_{r}=\{(t_{1},...,t_{n}):r>
t_{1}>...>t_{n}>0\},
$$
and introduce
$\bW^{m}_{r}$ as the closed linear
subspace of $L_{2}(\Omega,\cF,P)$ generated
by constants if $m=0$, by the set of
constants and
$$
\int_{0}^{r}f(t)\,dw^{k}_{t}
$$
if $m=1$
 or, if $m\geq 2$, by the set of constants and
$$
\int_{\Gamma^{n}_{r} }
f(t_{n},...,t_{1})
\,dw^{k_{n}}_{t_{n}}...\,dw^{k_{1}}_{t_{1}}
$$
$$
:=
\int_{0}^{r}dw^{k_{1}}_{t_{1}}\int_{0}^{t_{1}}dw^{k_{2}}_{t_{2}}
...\int_{0}^{t_{n-1}}f(t_{n},...,t_{1})\,
dw^{k_{n}}_{t_{n}},
$$
where $k$, $(k_{1},...,k_{n})$, $n\leq m$ are arbitrary and $f(t)$ and $f(t_{n},...,t_{1})$ are arbitrary Borel bounded functions of their arguments. The projection operator in $L_{2}
(\Omega,\cF,P)$ on $\bW^{n}_{r}$ we denote
by $\Pi^{n}_{r}$.

Consider the equation
\begin{equation}
                          \label{6.7.1}
 \cL u=0
\end{equation}
in $(-\infty,r)\times\bR^{d}$ 
with the terminal condition $u(r,x)=f(x)$,
where $f
\in C^{\infty}_{3}$. The results of
\cite{Kr_99} (where you throw away all stochastic terms) yield a solution of this
problem in certain class of functions $u$,
which, thanks to Corollary 4.12 of 
\cite{Kr_01_1},  are such that $u(t,\cdot)\in C^{\infty}_{3} $ for every $t$  
and each derivative in $x$ of any order is bounded uniformly on $[-T,r]\times \bR^{d}$
for any $T$. We denote
$$
T_{t,r}f(x)=u(t,x).
$$

Below we follow some arguments from
\cite{VK_76}.
Let $r>0$ and let $x_{s} $ be the solution of \eqref{6.15.2}
with $t=0$ and  $x=0$ . By It\^o's formula
applied to $T_{ s, r}f(x_{s})$ we get 
\begin{equation}
                                                     \label{6.16.3}
 f(x_{r})=   T_{0,r}f(0)+\int_{0}^{r}\sigma^{ik}D_{i}
   T_{ s,r}f(s,x_{s})\,dw^{k}_{s},
\end{equation}
where $\sigma^{ik}D_{i}
   T_{s,r}f(s,x )= \sigma^{ik}(s,x)D_{i}
   T_{s,r}f (x )$ and similar notation is also used below.

Next, we iterate \eqref{6.16.3}, that is we apply it to $\sigma^{ik}D_{i}
   T_{s,r}f(s,x )$ and $s$ in place of $f$
and $r$.
Then for any $s<r$ we get
$$
\sigma^{ik}D_{i}
   T_{s,r}f(s,x_{s})=T_{0,s}(\sigma^{ik}D_{i}
   T_{s,r}f)(0)
$$
\begin{equation}
                                                     \label{6.17.20}
+\int_{0}^{s}\sigma^{jm}D_{j}T_{u,s}\big(
\sigma^{ik}D_{i}
   T_{s,r}f\big)(u,x_{u})\,dw^{m}_{u}.
\end{equation}
After that we  substitute the result into \eqref{6.16.3}
to obtain
$$
f(x_{r})=T_{0,r }f(0) 
+\int_{0}^{r}T_{0,s}(\sigma^{ik}D_{i}
   T_{s,r}f)(0)\,dw^{k}_{s}
$$
\begin{equation}
                                                     \label{6.17.3}
+\int_{0}^{r}\Big(\int_{0}^{s}\sigma^{jm}D_{j}T_{u,s}\big(
\sigma^{ik}D_{i}
   T_{s,r}f\big)(u,x_{u})\,dw^{m}_{u}\Big)dw^{k}_{s}.
\end{equation}

Then we   repeat this procedure. Introduce  
\begin{equation}
                                                         \label{6.26.5}
Q^{k}_{s,r}f(x)=\sigma^{ik}(s,x)D_{i}T_{s,r}f(x).
\end{equation}
In this notation \eqref{6.16.3} and \eqref{6.17.3} become, respectively,
$$
 f(x_{r})=  
T_{0,r}f(0)+\int_{0}^{r}Q^{k_{1}}_{ t_{1},r}f(x_{t_{1}})\,dw^{k_{1}}_{t_{1}};
$$
$$
f(x_{r})=T_{0,r}f(0) 
+\int_{0}^{r}T_{0,t_{1}}Q^{k_{1}}_{ t_{1},r}f(0)\,dw^{k_{1}}_{t_{1}}
$$
$$
+\int_{0}^{r}\Big(\int_{0}^{t_{1}}Q^{k_{2}}_{t_{2},t_{1} } 
Q^{k_{1}}_{ t_{1},r}f (x_{t_{2}})\,dw^{k_{2}}_{t_{2}}\Big)dw^{k_{1}}_{t_{1}}.
$$

By induction we obtain that for any $n\geq1$   and $r\geq0$ ($t_{0}=r$)
$$
f(x_{r})=T_{0,r}f(0) 
+\sum_{m=1}^{n}\int_{\Gamma^{m}_{r}}T_{0,t_{m}}
Q^{k_{m}}_{t_{m},t_{m-1}}\cdot...\cdot
Q^{k_{1}}_{ t_{1},r}f(0)\,dw^{k_{m}}_{t_{m}}\cdot...\cdot dw^{k_{1}}_{t_{1}}
$$
\begin{equation}  
                         \label{6.18.1}
+\int_{\Gamma^{n+1}_{r}}
Q^{k_{n+1}}_{t_{n+1},t_{n}}\cdot...\cdot
Q^{k_{1}}_{ t_{1},r}f(x_{t_{n+1}})\,dw^{k_{n+1}}_{t_{n+1}}\cdot...\cdot dw^{k_{1}}_{t_{1}}.
\end{equation}

Obviously, the last term in \eqref{6.18.1}
is orthogonal to $\bW^{n}_{r}$. Therefore,
we have the following.

\begin{lemma}
                      \label{lemma 6.7.1}
For any $r>0$, $f\in C^{\infty}_{3} $,
$n\geq 1$,
and $\xi:=f(x_{r}(0,0))$
\begin{equation}
                               \label{6.7.2}
E_{0,0}|\xi-\Pi^{n}_{r}\xi|^{2}
$$
$$
=\int_{\Gamma^{n+1}_{r} }\sum_{k_{1},...,k_{n +1}}T_{0,t_{n+1}}
\big[
Q^{k_{n+1}}_{t_{n+1},t_{n}}\cdot...\cdot
Q^{k_{1}}_{ t_{1},r}f\big]^{2}(0 )\,d t_{n+1 } \cdot...\cdot d t_{1}.
\end{equation}
\end{lemma}

Next,   recall the notation
$f_{(\eta)}=\eta^{i}D_{i}f$  and consider the system consisting of
\eqref{6.15.2} and
\begin{equation}
                    \label{6.20.4}
\eta_{s}=\eta+\int_{0}^{s}\sigma^{k}_{(\eta_{s})}(t+r,x_{r})\,dw^{k}_{r}
+\int_{0}^{s}b_{(\eta_{r})}(t+r,x_{r})\,dr.
\end{equation}

We borrow some text from \cite{Kr_21}
with natural modifications caused by the fact that there the coefficients are independent of time and also we avoid using the somewhat mysterious function
$K_{0}(x)$.

As we know, \eqref{6.15.2} 
 has a unique solution which we denote by $x_{s}(t,x)$.
By substituting it into \eqref{6.20.4} we see that the coefficients
of \eqref{6.20.4} grow linearly in $\eta$ and hence 
\eqref{6.20.4} also has a unique solution which we denote by
$\eta_{s}(t,x,\eta)$. By the way, observe that equation \eqref{6.20.4}
is linear with respect to $\eta_{t}$. Therefore
$\eta_{s}(t,x,\eta)$ is a linear function of $\eta$.
For the uniformity of notation we set $x_{s}(t,x,\eta)=x_{s}(t,x)$.
It is also well known 
  (see, for instance, Sections 2.7 and 2.8 of
\cite{Kr_77}) that, as   functions of $x$ and $(x,\eta)$, the processes
$x_{t}(x)$ and $\eta_{t}(x,\eta)$
 are infinitely differentiable in an appropriate sense and
their derivatives satisfy the equations which are obtained by formal
differentiation of \eqref{6.15.2} and \eqref{6.20.4},
respectively.

As in (6.6) of \cite{Kr_21} we have 
that, if $f\in C^{\infty}_{0}$, then for any $s>0$
\begin{equation} 
                           \label{6.21.5}
\big(Ef(x_{s}(t,x))\big)_{(\eta)}
=E\Big(\big(f_{(\eta_{s}(t,x,\eta))}\big)(x_{s}(t,x))\Big).
\end{equation}
Here is the ``time dependent'' version of Lemma 6.1 of \cite{Kr_21} with 
 the same proof as there.

\begin{lemma}
                                                     \label{lemma 6.21.10}
Let $f(x,\eta)$ be infinitely differentiable and such that
each of its derivatives grows  as $|x|+|\eta|\to\infty$
not faster than polynomially. Let $T\in\bR$. Then for $t\leq T$, the function
$u(t,x,\eta):=Ef\big((x_{T-t},\eta_{T-t})(t,x,\eta)\big)$
is infinitely differentiable in $(x,\eta)$ and 
each of its derivatives is continuous in $t$ and is by absolute
value bounded on each finite   interval in
$(-\infty,T]$
by a constant times $(1+|x|+|\eta|)^{m}$
for some $m$. Furthermore, $u(t,x,\eta)$ is Lipschitz continuous with respect to $t$,  in $(0,T)\times \bR^{2d}$, (a.e.)
$\partial_{t}u(t,x,\eta)$ exists and
$$
0= \partial_{t}u(t,x,\eta)+ (1/2)\sigma^{ik}\sigma^{jk}(t,x)u_{x^{i}x^{j}} (t,x,\eta)
+\sigma^{ik}\sigma_{(\eta)}^{jk}(t,x)u_{x^{i}\eta^{j}} (t,x,\eta)
$$
$$
+(1/2)\sigma_{(\eta)}^{ik} \sigma_{(\eta)}^{jk}(t,x)u_{\eta^{i}\eta^{j}}(t,x,\eta)
+b^{i}(t,x)u_{x^{i}} (t,x,\eta)+b^{i}_{(\eta)}(t,x)u_{\eta^{i}} (t,x,\eta)
$$
\begin{equation}
                          \label{6.21.3}
=:\partial_{t}u(t,x,\eta)+\check \cL(t,x,\eta)u(t,x,\eta).
\end{equation}
\end{lemma}

Now comes an analogue of Lemma 6.3 of \cite{Kr_21} with a very similar proof,
whose origin is presented there.

\begin{lemma}
                    \label{lemma 6.21.01}
 Let $x,\eta\in\bR^{d}$, $r\in\bR$,  and let $f\in C^{\infty}_{0}$.
 Then for any $t<r $
$(t_{0}=r)$ 
$$
E\big[f_{(\eta_{r-t}(t,x,\eta))}(x_{r-t}(t,x))\big]^{2}
\geq\Big[(T_{t, r}f(x))_{(\eta)}\Big]^{2}
$$
\begin{equation} 
                         \label{6.21.6}
+\sum_{n=1}^{\infty}\sum_{k_{1},...,k_{n}}
\int_{\Gamma^{n}_{r-t}}\Big[\big(T_{t,t+t_{n}}Q^{k_{n}}_{t+t_{n},t+t_{n-1}}
\cdot...\cdot  Q^{k_{1}}_{t+t_{1}, r}f(x)\big)_{(\eta)}\Big]^{2}\,dt_{n}
\cdot...\cdot dt_{1}.
\end{equation}
\end{lemma}

Proof. For $t\leq r$ introduce the notation 
$$
\check T_{t,r}u(x,\eta)=Eu\big((x_{r-t},\eta_{r-t}) (t,x,\eta)\big).
$$ 
Then, similarly to \eqref{6.16.3}, 
by using Lemma \ref{lemma 6.21.10} and
applying It\^o's formula to
$\big(\check T_{ t+s,  r}u\big)\big((x_{s},\eta_{s})(t,x,\eta)\big)$
for smooth bounded
$u(x,\eta)$  by dropping for simplicity the arguments $t,x$ and $\eta$  
in $x_{\cdot}(t,x)$ and $\eta_{\cdot}(t,x,\eta)$, we get
$$
u(x_{r-t},\eta_{r-t})=\check T_{t, r}u(x,\eta)
$$
$$
+\int_{0}^{r-t}\Big[\sigma^{ik}(t+t_{1},x_{t_{1}})D_{x^{i}}
\check T_{t+t_{1}, r}u(x_{t_{1}},\eta_{t_{1}})$$
$$
+\sigma^{ik}_{(\eta_{t_{1}})}(t+t_{1},x_{t_{1}})D_{\eta^{i}}
\check T_{t+t_{1}, r}u(x _{t_{1} },\eta_{t_{1}} )\Big]\,dw^{k}_{t_{1}}.
$$
It follows that
$$
Eu^{2}(x_{r-t},\eta_{r-t})= \big(\check T_{t, r}u(x,\eta)\big)^{2}
$$
$$
+\sum_{k}\int_{0}^{r-t}E\Big[\sigma^{ik}(t+t_{1},x_{t_{1}})D_{x^{i}}
\check T_{t+t_{1}, r}u(x_{t_{1}},\eta_{t_{1}})
$$
\begin{equation}
                                                        \label{6.21.7}+
\sigma^{ik}_{(\eta_{t_{1}})}(t+t_{1},x_{t_{1}})D_{\eta^{i}}
\check T_{t+t_{1}, r}u(x _{t_{1}},\eta_{t_{1}})\Big]^{2}\,dt_{1}.
\end{equation}
By using Fatou's lemma, formulas like \eqref{6.21.5}, and well-known estimates of
the derivatives of solutions of It\^o's
equations with respect to initial data,
one easily carries \eqref{6.21.7},
with = replaced by $\geq$, over to smooth $u(x,\eta)$
whose derivatives have no more than polynomial growth
as $|x|+|\eta|\to\infty$. In particular, one can apply thus modified
\eqref{6.21.7} to $u(x,\eta)=f_{(\eta)}(x)$. Then,
after noting that, in light of \eqref{6.21.5}, in that case
$$
\sigma^{ik}(t+t_{1},x)D_{x^{i}}
\check T_{t+t_{1}, r}u(x,\eta )+
 \sigma^{ik}_{(\eta )}(t+t_{1},x)D_{\eta^{i}}
\check T_{t+t_{1}, r}u(x,\eta ) 
$$
$$
=\sigma^{ik}(t+t_{1},x )D_{x^{i}}(T_{t+t_{1}, r}f(x))_{(\eta)}
+\sigma^{ik}_{(\eta )}(x )D_{\eta^{i}}(T_{t+t_{1}, r}f(x))_{(\eta)}
$$
$$
=\big(\sigma^{ik}(t+t_{1},x )D_{x^{i}} T_{t+t_{1}, r}f(x)\big)_{(\eta)}
=\big(Q^{k}_{t+t_{1}, r}f(x)\big)_{(\eta)},
$$
we obtain
$$
E\big[f_{(\eta_{r-t} )}(x_{r-t} )\big]^{2}
\geq\Big[(T_{t, r}f(x))_{(\eta)}\Big]^{2}
+\sum_{k_{1}}
\int_{0}^{r-t}E\big[(Q^{k_{1}}_{t+t_{1}, r}f)_{(\eta_{t_{1}} )}
(x_{t_{1}} )\big]^{2}\,dt_{1}.
$$
By applying this formula to $Q^{k_{1}}_{t+t_{1}, r}f$ in place of $f$
we get   
$$
E\big[f_{(\eta_{t})}(x_{t})\big]^{2}
\geq\Big[(T_{t, r}f(x))_{(\eta)}\Big]^{2}+\sum_{k_{1}}
\int_{0}^{r-t} \big[(T_{t,t+t_{1}}Q ^{k_{1}}_{t+t_{1}, r}f (x  ))_{(\eta)}\big]^{2}\,dt_{1} 
$$
$$
+\sum_{k_{1},k_{2}}\int_{0}^{r-t}dt_{1}
\int_{0}^{t_{1}}E
\big[(Q^{k_{2}}_{t+t_{2},t+t_{1}}Q^{k_{1}}_{t+t_{1}, r}f)_{(\eta_{t_{2})}}(x_{t_{2}})
\big]^{2}\,dt_{2}.
$$

 Using induction shows that for any $n\geq1$ 
$$
E\big[f_{(\eta_{r-t} )}(x_{r-t} )\big]^{2}
\geq\Big[(T_{t, r}f(x))_{(\eta)}\Big]^{2}
$$
$$
+\sum_{m=1}^{n}\sum_{k_{1},...,k_{m}}
\int_{\Gamma^{m}_{r-t}}\Big[
I^{k_{1},...,k_{m}}(t_{1},...,t_{m})\Big]^{2}\,dt_{m}
\cdot...\cdot dt_{1}
$$
$$
+\sum_{k_{1},...,k_{n+1}}
\int_{\Gamma^{n+1}_{r-t}}E\Big[J^{k_{1},...,k_{n+1}}(t_{1},...,t_{n+1})\Big]^{2}\,dt_{n+1}
\cdot...\cdot dt_{1},
$$
where ($t_{0}=r$)
$$
I^{k_{1},...,k_{m}}(t_{1},...,t_{m})=
\big(T_{t,t+t_{m}}Q^{k_{m}}_{t+t_{m},t+t_{m-1}}
\cdot...\cdot  Q^{k_{1}}_{t+t_{1}, r}f(x)\big)_{(\eta)},
$$
$$
J^{k_{1},...,k_{n+1}}(t_{1},...,t_{n+1})=\big( Q^{k_{n+1}}_{t+t_{n+1 },t+t_{n}}
\cdot...\cdot  Q^{k_{1}}_{t+t_{1}, r}f\big)_{(\eta_{t_{n+1}})}
(x_{t_{n+1}})
$$
This yields \eqref{6.21.6} and proves the lemma.

In the proof of Theorem \ref{theorem 6.9.1} 
below the assumption
that $d\geq p_{D\sigma},  p_{b}>2$ is essential.  

\begin{theorem}
                          \label{theorem 6.9.1}
There exists $\beta=\beta(d,d_{1},\delta,   p_{D\sigma},  p_{b})>0$
such that if 
$$
\widehat{D\sigma}_{M},\hat b_{M}\leq \beta,
$$
 then
for any $f\in C^{\infty}_{0}$, $r>0$, $p\geq 2$, 
  $n>0$, and
\begin{equation}
                                 \label{6.11.5}
u(t,x,\eta)= E\big[f_{(\eta_{r}(t,x,\eta))}(x_{r}(t,x))\big]^{2},
\end{equation} 
we have  
$$
\int_{0}^{r} \Big( \int_{\bR^{2d}}
 u^{p}(t,x,\eta)  h(\eta)\,dxd\eta \Big)^{n}  dt
$$
\begin{equation}
                                 \label{6.9.1}
\leq N(d,\delta,p_{D\sigma},p_{b}, p,
\|D\sigma_{B}\|, \|b_{B}\|, r_{\sigma},r_{b},r
 ,n ) \Big( \int_{\bR^{d}}|Df|^{2p}\,dx
 \Big)^{n} ,
\end{equation}
where $h(\eta)=(1+|\eta|)^{-2p-d-1}$.

\end{theorem}

We know that $u(t,x,\eta)$ is a quadratic function with respect to $\eta$ and it satisfies
\eqref{6.21.3}, and we are going to prove
\eqref{6.9.1} by multiplying \eqref{6.21.3}
by $u^{p-1}(t,x,\eta)h  $ and integrating by parts. The possibility to do so is justified
by estimate (6.9) of \cite{Kr_21}:
for $t\in[0,r]$ and $x,\eta \in\bR^{d}$ 
$$
 |u(t,x,\eta)|+| u_{x}(t,x,\eta)|+
| u_{\eta}(t,x,\eta)|
$$
$$
+| u_{xx}(t,x,\eta)|
+| u_{x\eta}(t,x,\eta)| +| u_{\eta\eta}(t,x,\eta)| 
\leq M e^{-\mu |x| } (1+|\eta|^{2}) ,
$$
which holds for some constants $M, \mu>0$.

Observe that, for $s\in[0,r]$ and
$$
I(s):= \int_{\bR^{2d}} u^{p}(t,x,\eta)  h(\eta)\,dxd\eta  
$$
we have
$$
\frac{d}{ds}I^{n}(s)=nI^{n-1}
p\int_{s}^{r}\int_{\bR^{2d}}h(\eta)u^{p-1}(t,x,\eta)
\partial_{t}u(t,x,\eta)\,dxd\eta.
$$
It follows by Gronwall's inequality that,
to prove the theorem, it suffices to prove the following.

\begin{lemma}
                        \label{lemma 6.10.1}
There exists $\beta=\beta(d,d_{1},\delta,p_{D\sigma},  p_{b})>0$
such that if 
$$
\widehat{D\sigma}_{M},\hat b_{M}\leq \beta,
$$
 then
there is a constant $N= N(d,d_{1},\delta,p_{D\sigma},  p_{b}, r_{\sigma},p,r_{b})$ such that for any $t\in[0,r]$
$$
\int_{\bR^{2d}}h(\eta)u^{p-1}(t,x,\eta)
\check \cL(t,x,\eta)u(t,x,\eta)\,dxd\eta
$$
\begin{equation}
                                 \label{6.10.1}
   \leq N(1 +\widetilde{D\sigma}_{B}^{2}(t)
+\tilde b_{B}^{2}(t))\int_{\bR^{2d}}h(\eta)
u^{p }(t,x,\eta)\,dxd\eta.
\end{equation}

\end{lemma}

To prove Lemma \ref{lemma 6.10.1}, we use the following generalizations of
Lemma 4.1 of \cite{8}.

\begin{lemma}
                         \label{lemma 12.8.1}

For any $n\geq1$,  $p_{i} >0$, $i=1,...,n$,  and
$\kappa\geq d+1+p_{1}+...+p_{n}$ there exists a constant
$N=N(d,\kappa,n,p_{i})$ such that
for  any     $A_{i}\in \bR^{d}$, $i=1,...,n$, we have
\begin{equation}
                        \label{12.18.1}
|A_{1}|^{p_{1}}\cdot...\cdot|A_{n}|^{p_{n}}\leq 
N\int_{\bR^{d}}h(\eta)|(\eta,A_{1})|^{p_{1}}\cdot...\cdot |(\eta,A_{n})|^{p_{n}}
\,d\eta,
\end{equation}
where $h(\eta)=(1+|\eta|^{\kappa})^{-1}$
and $(\eta,A_{i})$ is the scalar product
of $\eta$ and $A_{i}$. Furthermore, if a $d\times d$-symmetric matrix $B$ is nonnegative
definite, $q\geq2$, and $h=(1+|\eta|)^{-2q-d-1},$ then
\begin{equation}
                        \label{6.9.2}
 \int_{\bR^{d}}h(\eta)|\eta|^{2}(\eta,B\eta)^{q-2}|B\eta|^{2}\,d\eta\leq
N(d,q)\int_{\bR^{d}}h(\eta)(\eta,B\eta)^{q}\,d\eta
\end{equation}
\end{lemma}

Here \eqref{12.18.1} is proved in \cite{8}.
Estimate \eqref{6.9.2} is true because its
left-hand side, obviously, less than
$N(d,q)\text{\rm tr}\, B$, and one estimates
$\text{\rm tr}\, B$ through the right-hand side
of \eqref{6.9.2} by dividing both by 
$\text{\rm tr}\, B$ and assuming that the estimate is wrong.

We also need a slight modification
of Lemma 3.5 of \cite{1}.  

\begin{lemma}
                         \label{lemma 8.10.2}
There exists   constants $\hat N=\hat N(d,p_{b})$ and $  N= N(d, p_{b},r_{b})$
such that
for any $w\in W^{1}_{2}(\bR^{d})$, $t\in\bR$,
and $\beta>0$
$$
\int_{\bR^{d}}|  b (t,x)|w(x)|Dw(x)|\,dx   
$$
$$
\leq
(\hat N \hat b_{M} +\beta)   
 \int_{\bR^{d}} |Dw  |^{2}\,dx +
(N \hat b_{M} +\beta^{-1}\tilde b _{B}^{2}(t))\int_{\bR^{d}} |w  |^{2}\,dx . 
$$
\end{lemma}  

Proof. Observe that $p_{b}>2$ 
 (Assumption \ref{assumption 6.3.10})  and by Lemma 3.5 of \cite{1} 
\begin{equation}
                            \label{6.11.1}
\int_{\bR^{d}} |b_{M}|^{2}|w|^{2}\,dx
\leq  \hat N \hat b_{M}^{2}  
 \int_{\bR^{d}} |Dw  |^{2}\,dx +
N \hat b_{M}^{2}\int_{\bR^{d}} |w  |^{2}\,dx.
\end{equation}
Hence,
$$
\int_{\bR^{d}}|  b _{M}|w|Dw|\,dx \leq
\Big(\int_{\bR^{d}} |Dw  |^{2}\,dx\Big)^{1/2}
\Big(\int_{\bR^{d}} |b_{M}|^{2}|w|^{2}\,dx\Big)^{1/2}    
$$
$$
\leq    \hat N \hat b_{M}   
 \int_{\bR^{d}} |Dw  |^{2}\,dx +
N \hat b_{M} \int_{\bR^{d}} |w  |^{2}\,dx .
$$
Upon combining this with
$$
\int_{\bR^{d}}|  b _{B}|w|Dw|\,dx  \leq
\tilde b _{B}\int_{\bR^{d}}|  w|\,|w_{x }|\,dx
\leq
\beta\int_{\bR^{d}} |w_{x} |^{2}\,dx+\beta^{-1}\tilde b _{B}^{2}
 \int_{\bR^{d}}   |w|^{2}\,dx ,   
$$
we get the result. The lemma is proved.

{\bf Proof of Lemma \ref{lemma 6.10.1}}.
We, actually, follow the proof of Lemma 7.4
of \cite{8} up to the places where we need to treat the terms  with   $b$ (and $D\sigma$)   which are now the sums of two and where
partial derivatives with respect to $\eta$
are involved.
In \cite{8} in a more general setting
to deal with partial derivatives with respect to $\eta$ an auxiliary  function
$K_{0}(x)$ was used. We start by 
  assuming  that 
$$
\widehat{D\sigma}_{M} \leq 
\beta,\quad \hat b_{M}\leq \beta
$$
  for some $\beta\in(0,1]$,  
which we are going to specify in the proof. 

For simplicity of notation we drop the arguments $t,x,\eta$. 
We also write $U\sim V$ if
their integrals over $\bR^{2d}$ coincide, and $U\prec V$ if the integral of
$U$ is less than or equal to that of $V$.  Below the constants called $N$  depend
 only on $d$, $d_{1}$, $\delta$,   $  p_{b},p,  p_{D\sigma}$, $ r_{\sigma}$,  $r_{b}$. The constants called $\hat N$
depend only on $d,d_{1},\delta  $,   $  p_{b}, p_{D\sigma}$. Define $h=(1+|\eta|)^{-2p-d-1}$.

Set   $w=u^{p/2}$ and note simple formulas:
$$
u^{p-1}u_{x}=(2/p)ww_{x},\quad u^{p-2}u_{x^{i}}u_{x^{j}}
=(4/p^{2})w_{x^{i}}w_{x^{j}}.
$$
Then
denote by $\check \cL_{1}$ the sum of the first-order terms in $\check \cL$
 and observe that integrating by parts shows that
$$
h u^{p-1}  b^{i}_{(\eta)} u_{\eta^{i}}  \sim -(1/p)
h _{  \eta^{i}}  b^{i}_{(\eta)}u^{p}
-(1/p)h b^{i}_{x^{i}} u^{r}
$$
$$
\sim   (2/p)\eta^{k}h _{\eta^{i} } b^{i}w  w_{x^{k}} 
+(2/p)h b^{i}ww_{x^{i}} .
$$
Hence,  
$$
hu^{p-1}\check \cL_{1}u \sim  
 (2/p)\eta^{k}h_{\eta^{i} } b^{i}w  w_{x^{k}} 
+(4/p)hb^{i}ww_{x^{i}} .
$$
Note that $|\eta|\,|h_{\eta}|\leq
N( d,p)h$, so that 
 by Lemma \ref{lemma 8.10.2} (recall $\hat b_{M}\leq \beta$)
\begin{equation}
                                                        \label{6.22.1}
 h u^{p-1}\check L_{1}u\prec  
 \hat N  \beta  h|w_{x} |^{2}+hN (1
+\beta^{-1}\tilde b_{B}^{2})|w|^{2}.
\end{equation}

Starting to deal with the second order derivatives notice that  
$$
h u^{p-1}(1/2)\sigma^{ik}\sigma^{jk} u_{x^{i}x^{j}} \sim -
((p-1)/2)u^{p-2}h\sigma^{ik}u_{x^{i}} \sigma^{jk} u_{x^{j} } 
$$
$$
-(1/2)h\big[\sigma^{ik}_{x^{i}}\sigma^{jk}+
\sigma^{ik}\sigma^{jk}_{x^{i}}\big]u^{p-1}u_{x^{j}} =
-((2p-2)/p^{2}) h\sigma^{ik}w_{x^{i}} \sigma^{jk} w_{x^{j} } 
$$
$$
-(1/p)h\big[\sigma^{ik}_{x^{i}}\sigma^{jk}+
\sigma^{ik}\sigma^{jk}_{x^{i}}\big]ww_{x^{j}}  
\leq -(1/d) h\sigma^{ik}w_{x^{i}} \sigma^{jk} w_{x^{j} }
$$
$$
+h\Big|\big[\sigma^{ik}_{x^{i}}\sigma^{jk}+
\sigma^{ik}\sigma^{jk}_{x^{i}}\big]ww_{x^{j}}\Big|,
$$
where the inequality (to simplify the writing) is due to the fact that $d\geq p\geq2$.
In this inequality the first term on the right is dominated in the sense
of $\prec$ by 
$$
-(1/d)\delta h|w_{x}|^{2}
$$
 (because $a\in\bS_{\delta}$).
The remaining term contains $ww_{x^{i}}$ and we treat it as above recalling that $\widehat{D \sigma}_{M} 
\leq\beta$. 
Then we get 
$$
h u^{p-1}(1/2)\sigma^{ik}\sigma^{jk} u_{x^{i}x^{j}}\prec
-\big[(1/d)\delta-\hat N \beta  \big] h|w_{x}|^{2}
$$
\begin{equation}
                                                        \label{6.22.3}
 +N(1+\beta^{-1}\widetilde{D\sigma}_{B}^{2})h|w|^{2}.
\end{equation}
 
To deal with the terms containing derivation
with respect to $\eta$ we prepare two inequalities.
Observe that as in \eqref{6.11.1}
with $\sigma^{k}$ denoting the $k$-th column of $\sigma$
\begin{equation}   
                                                          \label{7.5.2}
h| \sigma^{k}_{x}|^{2}w^{2}\,dx\prec \hat Nh\beta^{2}  |w_{x}|^{2} +Nh
 |w |^{2}  
\end{equation} 
and, owing to this and   Lemma  \ref{lemma 12.8.1},   
\begin{equation}
                             \label{6.11.2}
h|\eta|^{2} \sum_{k}|\sigma^{k}_{x}|^{2}
|w_{\eta}|^{2}\prec \hat N h  \sum_{k}|\sigma^{k}_{x}|^{2}
w ^{2}\prec
\hat N h\beta^{2}  
   |w_{x} |^{2} +
N h\beta^{2}  |w  |^{2}.
\end{equation}

Now,
$$
h u^{p-1}\sigma^{ik}\sigma_{(\eta)}^{jk} u_{x^{i}\eta^{j}}  \sim
-(p-1)h\sigma^{ik}u^{p-2} u_{\eta^{j}} \sigma_{(\eta)}^{jk} u_{x^{i}}
$$
$$
-u^{p-1} u_{x^{i}} \big[h_{\eta^{j}}  \sigma^{ik}\sigma_{(\eta)}^{jk}
+h \sigma^{ik}\sigma_{x^{j}}^{jk}]=-((4p-4)/p^{2})
h\sigma^{ik}  w_{\eta^{j}} \sigma_{(\eta)}^{jk} w_{x^{i}}
$$
$$
-(2/p)ww_{x^{i}} \big[h_{\eta^{j}}  \sigma^{ik}\sigma_{(\eta)}^{jk}
+h\sigma^{ik}\sigma_{x^{j}}^{jk}].
$$
We estimate the first term on the right roughly using
$$
|\sigma^{ik} w_{\eta^{j}} \sigma_{(\eta)}^{jk} w_{x^{i}}|\leq
\beta|w_{x}|^{2}+\hat N\beta^{-1} |\eta|^{2} \sum_{k}|\sigma^{k}_{x}|^{2}
|w_{\eta}|^{2}
$$
and then applying \eqref{6.11.2}.
The second term  contains $ww_{x^{i}}$ and allows the same handling as before, using Lemma \ref{lemma 8.10.2}.
 Therefore,
\begin{equation}   
                        \label{6.22.4}
h u^{p-1}\sigma^{ik}\sigma_{(\eta)}^{jk} u_{x^{i}\eta^{j}} \prec
   \hat N \beta    h|w_{x}|^{2}
+N(1+\beta^{-1}\widetilde{D\sigma}_{B}^{2} )h|w|^{2}  .
\end{equation}

The last term in $h u^{p-1}\check  Lv$ containing $\sigma$ is  
$$
h u^{p-1}(1/2)\sigma_{(\eta)}^{ik} \sigma_{(\eta)}^{jk} u_{\eta^{i}\eta^{j}} \sim
-((p-1)/2)h\sigma_{(\eta)}^{ik} u^{p-2} u_{ \eta^{j}}  \sigma_{(\eta)}^{jk}u_{\eta^{i}}
$$
$$
-(1/2)u^{p-1}\sigma_{(\eta)}^{ik} u_{\eta^{i} }\big[h_{\eta^{j}}  \sigma_{(\eta)}^{jk}
+h \sigma_{x^{j}}^{jk}\big]-(1/(2p))h (u^{p})_{\eta^{i}}
 \sigma_{x^{j}}^{ik} \sigma_{(\eta)}^{jk}
$$
$$
\prec \hat NJ +I,
$$
where  
$$
J=h(|\eta|^{2}|w_{\eta}|^{2}+w^{2})\sum_{k}|\sigma^{k}_{x}|^{2},
$$
which, according to \eqref{6.11.2},
is dominated in the sense of $\prec$ by 
$
\hat N h\beta^{2} 
   |w_{x} |^{2} +
N h\beta^{2}  |w  |^{2}
$,
and
$$
I=-(1/(2p))h(w^{2})_{\eta^{i}}
 \sigma_{x^{j}}^{ik} \sigma_{(\eta)}^{jk}
\sim
(1/(2p))w^{ 2}\sigma_{x^{j}}^{ik}\big[h_{\eta^{i}} \sigma_{(\eta)}^{jk}
+h\sigma_{x^{i}}^{jk}\big]
$$
$$
\prec \hat Nh\sum_{k}|\sigma^{k}_{x}|^{2}w^{2}
\prec \hat N h\beta^{2} 
   |w_{x} |^{2} +
N h\beta^{2}  |w  |^{2}.
$$

Hence, 
\begin{equation}
                                                                 \label{7.9.3}
h u^{p-1}(1/2)\sigma_{(\eta)}^{ik} \sigma_{(\eta)}^{jk} u_{\eta^{i}\eta^{j}}
\prec   \hat Nh\beta^{2}   |w_{x}|^{2}+Nh
\beta^{2} |w|^{2}.
\end{equation}

By combining \eqref{6.22.1}, \eqref{6.22.3},
 \eqref{6.22.4}, and \eqref{7.9.3},  and taking into account that $\beta\leq 1$,
 we see that  
 $$
h u^{q-1}\check Lv\prec 
-\big[(1/d)\delta-\hat N \beta 
 \big] h|w_{x}|^{2} +N(1 
+\beta^{-1}\widetilde{D\sigma}_{B}^{2}
+\beta^{-1}\tilde b_{B}^{2})h|w|^{2}.
$$
Upon choosing $\beta\leq \delta/(\hat Nd)$
we get \eqref{6.10.1}. The lemma is proved.
\qed

Now we   state the key result on our way of
  justifying
 the Conjecture.
One should pay attention to what $N$ in
\eqref{6.11.4} depends on.

  Recall that Assumption \ref{assumption 6.3.10} is supposed to be satisfied in which we,
in particular,  assume that $ p(d,p_{b})
\geq 2, p(d,p_{b})>d/2$.  
\begin{theorem}
                    \label{theorem 6.11.1}  
Take $\beta$ from Theorem \ref{theorem 6.9.1}
and assume that $\widehat{D\sigma}_{M},\hat b_{M}\leq\beta$.  Take $f\in C^{\infty}_{0}$
and $r>0$  and set $x_{r} =x_{r}(0,0)$, $\xi=f(x_{r} )$. Then for  $p =p(d,p_{b})$   
 we have
\begin{equation}
                            \label{6.11.4}
\sum_{n=1}^{\infty}E|\xi-\Pi^{n}_{r}\xi|^{2}
\leq N(d,\delta,p_{D\sigma},p_{b},p, \| D\sigma_{B}\|,
 \| b_{B}\|, r_{\sigma},r_{b},r)\Big(\int_{\bR^{d}}|Df|^{2p }\,dx\Big)^{1/p}.
\end{equation}
\end{theorem}

Proof.   Since $u$ in \eqref{6.11.5} is quadratic with respect to $\eta$, estimate
\eqref{6.9.1} implies that  
$$
\int_{0}^{r} \Big( \int_{\bR^{d}}\sup_{|\eta|\leq 1}
u^{p}(t,x,\eta)\,dx \Big)^{n} dt\leq N \Big( \int_{\bR^{d}}|Df|^{2p}\,dx
 \Big)^{n} ,
$$
 where $n=n(p)$ is any number such that
$n>2/(2p-d)$, 
which, in turn, implies that
\begin{equation}
                                  \label{3.4.4}
\int_{0}^{r} \Big( \int_{\bR^{d}}g^{p}(s,x)\,d x
 \Big)^{q/p} ds\leq N \Big( \int_{\bR^{d}}|Df|^{2p}\,dx   \Big)^{q/p} ,
\end{equation}
where  $q=np$ and 
$$
g (s,x)=
\sum_{k}u\big(s,x,\sigma^{k}(s,x)\big).
$$
Since  $d/p+2/q<2$ and $p<p_{b}$   by Theorem \ref{theorem 6.4.2} 
 we obtain  
$$
\int_{0}^{r}T_{0,s}g(s,\cdot)(x)\,ds=
E\int_{0}^{r}g(s,x_{s})\,ds
\leq N\Big(\int_{\bR^{d}}|Df|^{2p}\,dx\Big)^{1/p}.
$$
After that it only remains to observe that,
in light of Lemma 3.4, we have
$$
g (s,x)\geq 
\sum_{n=1}^{\infty}\sum_{k_{1},...,k_{n+1}}
\int_{\Gamma^{n}_{r-s}}\Big[ Q^{k_{n+1}}_{s,s+ t_{n}}Q^{k_{n}}_{ s+t_{n},s+ t_{n-1}}
\cdot...\cdot  Q^{k_{1}}_{ s+t_{1}, r}f(x) \Big]^{2}\,dt_{n}
\cdot...\cdot dt_{1},
$$
so that
$$
\int_{0}^{r}T_{0,s}g(s,\cdot)(x)\,ds
$$
is greater than or equal to
$$
\sum_{n=1}^{\infty}\sum_{k_{1},...,k_{n+1}}\int_{0}^{r}
\int_{\Gamma^{n}_{r-s}}T_{0,s}\Big[ Q^{k_{n+1}}_{s,s+ t_{n}}Q^{k_{n}}_{ s+t_{n},s+ t_{n-1}}
\cdot...\cdot  Q^{k_{1}}_{ s+t_{1}, r}f(x) \Big]^{2}\,dt_{n}
\cdot...\cdot dt_{1}ds
$$
and the latter coincides with the left-hand side of \eqref{6.11.4} by Lemma
\ref{lemma 6.7.1}. The theorem is proved.\qed

 The hope is that if
we approximate $\sigma$ and $b$ with $\sigma^{m}$
 and $b^{m}$
for which \eqref{6.11.4} holds and call $x^{m}_{t}$
the corresponding solution, then in some sense
$x^{m}_{t}$ will converge to $x^{0}_{t}$ which solves
the original equation and for $\xi^{m}=f(x^{m}_{t})$
we will have $\xi^{m}-\Pi^{n}_{t}\xi^{m}\to
\xi^{0}-\Pi^{n}_{t}\xi^{0}$, which will allow
to prove \eqref{6.11.4} for thus obtained $x^{0}_{t}$.
Since \eqref{6.11.4} implies that the 
$\cF^{w}_{t}$-measurable variables $\Pi^{n}_{t}\xi$
converge to $\xi$, we would be able to conclude that
$x^{0}_{t}$ is a strong solution. Despite
the fact that the author failed a few times to 
rigorously justify the above program with different
approaches, there is
still hope that it is realistic.

Approximating $\sigma$ with smooth ones presents
some difficulties. 
Indeed, if $d=1$ and $\sigma(x)=\text{sign}\, x$,
then any mollification of $\sigma$ will have
zeros and the corresponding diffusion is degenerate.
Under our ``smoothness'' assumptions this never happens.

We need an additional restriction on $\widehat{D\sigma}_{M}$.
Take   nonnegative $\xi\in C^{\infty}_{0}(\bR)$, $\eta\in C^{\infty}_{0}(\bR^{d})$
with unit integrals and supports in unit balls, set $\zeta(t,x)=\xi (t)*\eta (x)$, $\zeta_{n}(t,x)=n^{d+1}\zeta(nt,nx)$,
$\eta_{n}( x)=n^{d }\eta( nx)$ and  define
 $$
\sigma^{(n)}=\sigma * \eta_{n}, 
\quad D\sigma^{(n)}_{M}=(D\sigma_{M}) * \eta_{n}, \quad
D\sigma^{(n)}_{B}=(D\sigma_{B}) * \eta_{n}, 
$$
where the convolution is performed with respect to $x$, and 
$$
b^{(n)}=b * \zeta_{n},\quad
b^{(n)}_{M}=b_{M} * \zeta_{n},\quad
b^{(n)}_{B}=b_{B} * \zeta_{n},
$$
where the convolution is performed with respect to $(t,x)$.
As is easy to see, for each $n$, $\sigma^{(n)}$
is bounded along with its any derivative of any order with respect to $x$, $b^{(n)}$
 is a smooth bounded function.
Furthermore, Minkowski's inequality easily
shows that $D\sigma^{(n)},b^{(n)}$ (with the same  $p_{D\sigma},p_{b}$) satisfy
$$
\widehat{D\sigma}^{(n)}_{M}\leq \widehat{D\sigma}_{M},\quad 
\widetilde{D\sigma^{(n)}}_{B}
 \leq 
\widetilde{D\sigma} _{B},
$$
$$
\hat b^{(n)}_{M}
\leq \hat b_{M}, \quad \beta_{b^{(n)}}
\leq \beta_{b},\quad \| b^{(n)}_{B}\|\leq \| b_{B}\|.
$$
Finally, introduce
$$
\Gamma_{m}=\{t: \widetilde{D\sigma} _{B}(t)\leq
m\},
$$
$$
\sigma^{(n)}_{m}(t,x)=\sigma^{(n)} (t,x)
I_{\Gamma_{m}}(t)+\kappa I_{\Gamma^{c}_{m}}(t),
$$
where $\kappa$ is any fixed $d\times d_{1}$-matrix such that $\kappa \kappa^{*}=(\delta^{ij})$.

\begin{lemma}
                    \label{lemma 6.13.1}
Set $a^{n}_{m}=\sigma^{(n)}_{m}\sigma^{(n)*}_{m}$. Then there is a sequence $m(n)\to\infty$
as $n\to\infty$ such that
for sufficiently large $n$ the eigenvalues of
$a^{n}_{m(n)}$ are between $\delta/4$ and $4/\delta$
if 
\begin{equation}
                              \label{6.14.2}
N(d )\widehat{D\sigma}_{M}\leq 
\delta^{1/2}/4,
\end{equation}
 where $N(d )$ is specified in the proof.
\end{lemma}

Proof. Obviously, we only need to deal with $t\in
\Gamma_{m(n)}$.
Observe that for such $t$
$$
|\sigma^{(n)*}_{m} \lambda|=|\sigma^{(n)*} \lambda|\leq \eta_{ n}
  * |\sigma^{*} \lambda|\leq\delta^{-1/2}|\lambda|.
$$
Therefore, we need only prove that for sufficiently large $n$ and $t\in \Gamma_{m(n)}$
\begin{equation}
                                                    \label{7.2.1}
|\sigma^{(n)*}(t,x)\lambda|\geq |\lambda|\delta^{ 1/2}/ 2.
\end{equation}   
For any $ y$ we have  
$$
|\sigma^{(n)*}(t,x)\lambda|\geq |\sigma^{*}(t ,x-y)\lambda|-
 |(\sigma^{(n)*}(t,x)-\sigma^{*}(t ,x-y))\lambda|
$$
$$
\geq |\lambda|\big(\delta^{1/2}
-|\sigma^{(n)*}(t,x)-\sigma^{*}(t ,x-y)|\big).
$$
Furthermore, by using Poincar\'e's inequality
and recalling that $t\in \Gamma_{m(n)}$, we obtain
$$
I:=\int_{\bR^{d }}| \sigma^{(n)*} (t,x)-\sigma^{*}(t,x-y)  |\eta_{n}(y)\,dyds
$$
$$
\leq\int_{B_{1}}\int_{B_{1}}
| \sigma (t, x-z/n)-\sigma (x-y/n)  |\eta (y)\eta (z)\,dydz
$$
$$
\leq N(d)\dashint_{B_{1/n}(x)}\dashint_{B_{1/n}(x)}|\sigma(t,y)-\sigma(t,z)|\,dydz
$$
$$
\leq N(d )n^{-1}\dashint_{B_{1/n}(x)}|D\sigma(t,y)|\,dy=: N(d )I_{1}+N(d )I_{2},
$$
where
$$
I_{1}=n^{-1}\dashint_{B_{1/n}(x)}|D\sigma_{M}(t,y)|\,dy
$$
$$
\leq n^{-1}\Big(\dashint_{B_{1/n}(x)}|D\sigma(t,y)|^{p_{D\sigma}}\,dy\Big)^{1/p_{D\sigma}}\leq \widehat{D\sigma}_{M}
$$
if $1/n\leq  r_{\sigma}$, and
$$
I_{2}=n^{-1}\dashint_{B_{1/n}(x)}|D\sigma_{B}(t,y)|\,dy\leq n^{-1}m(n).
$$
Hence,
$$
I\leq N(d )\widehat{D\sigma}_{M}+
N(d)n^{-1}m(n) 
$$
and we get the desired result by defining
$m(n)$ so that $N(d)n^{-1}m(n)\leq \delta^{1/2}/4$.
 The  lemma is proved.

 {\bf Acknowledgments}. The author is grateful
 to the handling editor who managed 
 in thirteen months to get referees' reports
  that helped improve the presentation.

\end{document}